\documentclass[12pt]{amsart}

\newtheorem{theorem}{Theorem}[section]

\newtheorem{corollary}[theorem]{Corollary}

\newtheorem{lemma}[theorem]{Lemma}

\newtheorem{proposition}[theorem]{Proposition}

\theoremstyle{definition}

\newtheorem{definition}[theorem]{Definition}

\newtheorem{remark}[theorem]{Remark}

\theoremstyle{parrafo}

\newtheorem{parrafo}[theorem]{{\!}}

\numberwithin{equation}{theorem}

\newcommand{\vMax}{{\underline{\mbox{Max} }  } \ }

\newcommand{\vword}{{\mbox{w-ord}}}

\newcommand{\nat}{\mathbb N}

\newcommand{\calo}{{\mathcal {O}}}

\newcommand{\Sing}{\mbox{Sing\ }}
\newcommand{\ord}{\mbox{ord}}

\title[]{ Elimination with applications to singularities
in positive characteristic.}

\author{Orlando Villamayor U.}

\address{Dpto. Matem\'aticas, Facultad de Ciencias, Universidad
Aut\'onoma de Madrid, Canto Blanco 28049 Madrid, Spain.}
\email{villamayor@uam.es}
\thanks{2000 {\em Mathematics subject classification. 14E15.}}
\thanks{.}

\subjclass{}

\keywords{Integral closure. Rees algebras.} \date{May 2006}
\dedicatory{To Professor H. Hironaka} \commby{}

\begin{document}
\maketitle
\begin{abstract} We present applications of elimination theory to the study
of singularities over arbitrary fields. A partial extension of a
function, defining resolution of singularities over fields of
characteristic zero, is discussed here in positive characteristic.

\end{abstract}

\tableofcontents






\part{ Introduction.}

\label{introduction}

Hironaka's theorem of embedded desingularization was proven by
induction on the dimension of the ambient space. This form of
induction is based on a reformulation of the resolution problem,
as a new resolution problem, but now in a smooth hypersurface of
the ambient space. Smooth hypersurfaces playing this inductive
role are called hypersurfaces of maximal contact. In the case of
resolution of embedded schemes defined by one equation,
hypersurfaces of maximal contact can be selected via a
Tschirnhausen transformation of the equation. However this
strategy for induction on resolution problems holds exclusively
over fields of characteristic zero, and fails over fields of
positive characteristic.

The objective of this paper is to discuss
results that grow from a different approach to induction, based on a
a form of elimination which holds over fields of arbitrary
characteristic.

Over fields of characteristic zero Hironaka proves that resolution
of singularities is achieved by blowing up, successively, at
smooth centers. Constructive resolution of singularities is a form
of desingularization where the centers are defined by an
upper semi-continuous function. The singular locus is stratified
by the level sets of the function. The closed stratum,
corresponding to the biggest value achieved by the function, is
the smooth center to be blown-up. Then a new upper semi-continuous
function is defined at the blow-up, which, in the same way,
indicates the next center to blow-up; and so on.

In this paper we show that there is a partial extension to arbitrary characteristic of the upper semi-continuous function in
\cite{Villa89}, defined there over fields of characteristic
zero (see Theorem  \ref{thbv}  and Proposition  \ref{prop617} ).
The notion of eliminations algebras, introduced in \cite{VV4}, will be used as a substitute for the notion of maximal contact. A second ingredient for this extension is Hironaka's Finite Presentation Theorem (p.119, \cite{Hironaka05}) (see \ref{hfpt}).

This partial extension of the function to positive characteristic
provides, in a canonical manner, a procedure of transformation of
singularities into singularities of a specific simplified form (with
``monomial'' elimination algebra). Over fields of characteristic
zero this is the well know reduction to the {\em monomial case} (see
\ref{llmonn}).



Hironaka defines a class of objects (couples), consisting of an ideal and a
positive integer. On this class he introduces two notions of equivalence.
The first equivalence is defined in terms of integral closure of ideals, and an equivalence class is called an {\em idealistic exponent}.


In Sections 1 and 2 we
give an overview of the main results in \cite{VV1}, where
idealistic exponents are expressed as Rees algebras, and where
this notion of equivalence of couples is reinterpreted in terms of
integral closure of Rees algebras.

In Section 3 we discuss Rees algebras with an
action of differential operators (Diff-algebras). We also reformulate Giraud's Lemma of differential
operators and monoidal transformations, in terms of Rees-algebras.

In Section 4 we recall the main ingredients that appear in the
definition of the upper-semi-continuous  stratifying function
mentioned above, and show that there is a very natural extension of
these functions to the class of Rees algebras.

The second notion of equivalence, called weak equivalence, is discussed here in Section 5); together with the Finite Presentation Theorem, which is a bridge among both notions of equivalence. Weak equivalence played a central role in definition of the stratifying upper semi-continuous function over fields of characteristic zero, and in proving the properties studied in \cite{Villa92}. Namely, the compatibility of constructive resolution with \'etale topology, smooth maps, and the property of equivariance.

The partial extension of this stratifying function to positive characteristic, which we finally address in Section 6), makes use of Hironaka's Finite Presentation Theorem, together with elimination of Diff-algebras as a substitute for maximal contact.

It is the context of Diff-algebras where our form of elimination is defined, and Diff-algebras are Rees-algebras enriched with the action of higher differential operators. Rees algebras extend to a Diff-algebras, and this extension is
naturally compatible with integral closure of algebras (\cite{VV1}). This interplay of Diff-algebras and integral closure is studied by
Kawanoue in \cite{kaw}, and in \cite{MK}, papers which present new ideas and technics in positive characteristic, and also provide an
upper semi-continuous function with a different approach.

We refer to \cite{Hironaka06}  for a program of Hironaka for
embedded resolution over fields of positive characteristic. We also
refer to \cite{Cut} and \cite{CP} for new proofs on non-embedded
resolution of singularities of schemes of dimension 3 and positive
characteristic.

I am grateful to the referee for many useful suggestions and
the careful reading of the paper. I profited from discussions with Ang\'elica Benito, Ana Bravo, Mari Luz Garc\'ia, and Santiago Encinas.

\section{Idealistic exponents and Rees algebras.}\label{sec1}

\begin{parrafo}\label{parsec}
{\em In what follows $V$ denotes a smooth scheme over a field $k$. A
{\em couple} $(J,b)$ is a pair where $J$ non-zero sheaf of ideals in
$\calo_V$, and $b$ is a positive integer. We will consider the class of all couples,
and transformation among them.

$\bullet$ Given a couple $(J,b)$,   the  {\bf closed set}, or {\bf singular locus},  is:
$$\Sing(J,b)=\{ x\in V / \nu_x(J_x)\geq b \},$$ namely the set of
points in $V$ where $J$ has order at least $b$ (here $\nu_x$
denotes the order at the local regular ring $ \calo_{V,x}$). The
set $\Sing(J,b)$ is closed in $V$.

$\bullet${\bf Transformation} of $(J,b)$:

 Let $Y \subset \Sing(J,b)$ be a closed and smooth subscheme, and let
\begin{equation*}
\begin{array}{ccccc}
 & V & \stackrel{\pi}{\longleftarrow} & V_{1} \supset
H=\pi^{-1}(Y)& \\ & Y & &  &\\
\end{array}
\end{equation*}
denote the monoidal transformation at $Y$. Since $Y \subset
\Sing(J,b)$ the total transform, say  $J\calo_{V_1}$,  can be
expressed as a product: $$J\calo_{V_1}= I(H)^b J_1$$ for a
uniquely defined $J_1 $ in $ \calo_{V_1}$. The new couple $(J_1, b)$
is called the {\em transform} of $(J, b)$. We denote the
transformation by:
\begin{equation}\label{ecnwc}
\begin{array}{ccccc}
 & V & \stackrel{\pi}{\longleftarrow} & V_{1} ,& \\ & (J,b) & & (J_1,b) &\\
\end{array}
\end{equation}
and a sequence of transformations by:

\begin{equation}\label{sectransp}
\begin{array}{cccccccc}
 & V & \stackrel{\pi_1}{\longleftarrow} & V_{1} & \stackrel{\pi_2}{\longleftarrow}
 &\ldots& \stackrel{\pi_k}{\longleftarrow}& V_k.\\ & (J,b) & & (J_1,b) &&&&(J_k,b)\\
\end{array}
\end{equation}
Let $H_i$ denote the exceptional hypersurface introduced by $\pi_i$, $1\leq i \leq k$, which we also consider as hypersurfaces in $V_k$ (by taking strict transforms). Note that in such case
\begin{equation}\label{express}
J \calo_{V_k }= I(H_1)^{c_1}\cdot I(H_2)^{c_2} \cdots
I(H_k)^{c_k}\cdot J_k
\end{equation}
for suitable exponents $c_2,\dots,c_k$, and $c_1=b$. Furthermore,
all $c_i = b$ if for every index $ i<k$ the center $Y_i $ is not
included in $ \cup_{j \leq i}H_j \subset V_i$ (the exceptional
locus of $V \longleftarrow V_i$).
The previous sequence is said to be a {\em resolution} of $(J,b)$ if:

1) $\Sing(J_k,b)=\emptyset$, and

2) $ \cup_{j \leq k}H_j \subset V_k$ has normal crossings.

So if (\ref{sectransp}) is a resolution, then $J_k$ has at most
order $b-1$ at points of $V_k$.

Of particular importance for resolution of singularities is the
case in which $J_k$ has order at most zero, namely when
$J_k=\calo_{V_k}$. In such case we say that (\ref{sectransp}) is a
{\em Log-principalization} of $J$.

Given \((J_{1},b_{1})\) and \((J_{2},b_{2})\), then
\begin{equation*}
\Sing(J_{1},b_{1})\cap\Sing(J_{2},b_{2})=\Sing(K,c)
\end{equation*}
where \(K=J_{1}^{b_{2}}+J_{2}^{b_{1}}\), and \(c=b_{1}\cdot
b_{2}\). Set formally \((J_{1},b_{1})\odot (J_{2},b_{2})=(K,c)\).

If \(\pi\) is permissible for both \((J_{1},b_{1})\) and
\((J_{2},b_{2})\), then it is permissible for \((K,c)\).
Moreover, if \((J'_1,b_{1})\), \((J'_2,b_{2})\), and \((K',c)\)
denote the transforms, then \((J'_1,b_{1})\odot
(J'_2,b_{2})=(K',c)\).

}
\end{parrafo}

\begin{parrafo}\label{paris}

{\em We now define a {\em Rees algebra} over $V$ to be a graded
noetherian subring of $\calo_V[W]$, say:
$$\mathcal{G}=\bigoplus_{k\geq 0}I_kW^k,$$
where $I_0=\calo_V$ and each $I_k$ is a sheaf of ideals. We assume
that at every affine open $U(\subset V)$, there is a finite set
$\mathcal{F}=\{f_1W^{n_1},\dots ,f_sW^{n_s}\},$ $n_i\geq 1$ and
$f_i\in \calo_V(U)$, so that the restriction of $\mathcal{G}$ to
$U$ is
$\calo_V(U)[f_1W^{n_1},\dots ,f_sW^{n_s}] (\subset
\calo_V(U)[W]).$

To a Rees algebra $\mathcal{G}$ we attach a closed set:
$$\Sing(\mathcal{G}):=\{ x\in V/ \nu_x(I_k)\geq k, \mbox{ for each }
k\geq 1\},$$ where $\nu_x(I_k)$ denotes the order of the ideal
$I_k$ at the local regular ring $\calo_{V,x}$.
}
\end{parrafo}

\begin{remark}\label{rkK1} Rees algebras are related to Rees rings. A Rees algebra is a Rees ring if, given an affine open set
$U\subset V$, $\mathcal{F}=\{f_1W^{n_1},\dots ,f_sW^{n_s}\}$
can be chosen with all degrees $n_i=1$.
 Rees algebras are integral closures of Rees rings in a
suitable sense. In fact, if $N$ is a positive integer divisible by
all $n_i$, it is easy to check that $$\calo_V(U)[f_1W^{n_1},\dots
,f_sW^{n_s}]=\oplus_{r\geq 0}I_rW^r (\subset \calo_V(U)[W]),$$ is
integral over the Rees sub-ring $\calo_V(U)[I_NW^N](\subset
\calo_V(U)[W^N])$.

\end{remark}

\begin{proposition}\label{prop1} Given an affine open $U\subset V$, and
$\mathcal{F}=\{f_1W^{n_1},\dots ,f_sW^{n_s}\}$ as above,
$$\Sing(\mathcal{G})\cap U= \cap_{1\leq i \leq s}\{ord (f_i) \geq n_i\}.$$

\end{proposition}
\proof
 Since $\nu_x(f_i)\geq n_i$ for $x\in
\Sing(\mathcal{G})$, $0\leq i \leq s$;
$$\Sing(\mathcal{G})\cap
U\subset  \cap_{1\leq i \leq s}\{ord (f_i) \geq n_i\}.$$

On the other hand, for any index $N\geq 1$, $I_N(U)W^N$ is
generated by elements of the form $G_N(f_1W^{n_1},\dots
,f_sW^{n_s})$, where $G_N(Y_1,\dots ,Y_s)\in \calo_U[Y_1,\dots
,Y_s]$ is weighted homogeneous of degree $N$, provided each $Y_j$
has weight $n_j$. The reverse inclusion is now clear.
\endproof

\begin{parrafo}\label{14S}{\em
A monoidal transformation of $V$ on a smooth sub-scheme $Y$, say $V\stackrel{\pi}{\longleftarrow}
V_1$ is said to be {\em permissible} for
$\mathcal{G}$ if $Y\subset \Sing(\mathcal{G})$. In such case, for
each index $k\geq 1$, there is a sheaf of ideals, say
$I_k^{(1)}\subset \calo_{V_1}$, so that
$I_k\calo_{V_1}= I(H)^k I_k^{(1)},$ where $H$ denotes the exceptional locus of $\pi$.
One can easily check that
$$\mathcal{G}_1=\bigoplus_{k\geq 0}I^{(1)}_kW^k$$
is a Rees algebra over $V_1$, which we call the {\em transform} of
$\mathcal{G}$, and denote by:

\begin{equation}
\begin{array}{ccccc}
 & V & \stackrel{\pi}{\longleftarrow} & V_{1} & \\ & \mathcal{G} & & \mathcal{G}_1 &\\
\end{array}
\end{equation}

A sequence of transformations will be denoted as

\begin{equation}\label{sectransp11}
\begin{array}{cccccccc}
 & V & \stackrel{\pi_1}{\longleftarrow} & V_{1} & \stackrel{\pi_2}{\longleftarrow}
 &\ldots& \stackrel{\pi_k}{\longleftarrow}&V_k.\\ & \mathcal{G} & & \mathcal{G}_1 &&&&\mathcal{G}_k\\
\end{array}
\end{equation}
}
\end{parrafo}

\begin{definition}\label{defresg}

Sequence (\ref{sectransp11}) is said to be a {\em resolution}
of $\mathcal{G}$ if:

1) $\Sing(\mathcal{G}_k)=\emptyset$.

2)The union of the exceptional components, say $ \cup_{j \leq k}H_j \subset V_k$, has normal crossings.

\end{definition}
\begin{parrafo}\label{amalg}{\em

Given two Rees algebras over \(V\), say
\(\mathcal{G}_{1}=\bigoplus_{n\geq 0}I_{n}W^n\) and
\(\mathcal{G}_{2}=\bigoplus_{n\geq 0}J_{n}W^n\), set
\(K_n=I_n+J_n\) in \(\calo_V\), and define:
\begin{equation*}
\mathcal{G}_{1}\odot\mathcal{G}_{2}=\bigoplus_{n\geq 0}K_{n}W^n,
\end{equation*}
as the subalgebra of \(\calo_V[W]\) generated by \(\{K_{n}W^n,
n\geq 0\}\).

Let $U$ be an affine open set in $V$. If the restriction of
$\mathcal{G}_1$ to $U$ is $\calo_V(U)[f_1W^{n_1},\dots
,f_sW^{n_s}]$, and that of $\mathcal{G}_2$ is
$\calo_V(U)[f_{s+1}W^{n_{s+1}},\dots ,f_tW^{n_t}]$, then the
restriction of $\mathcal{G}_{1}\odot\mathcal{G}_{2}$ is
$$\calo_V(U)[f_1W^{n_1},\dots ,f_sW^{n_s}, f_{s+1}W^{n_{s+1}},\dots
,f_tW^{n_t}].$$

One can check that:
\begin{enumerate}
    \item \(\Sing(\mathcal{G}_{1}\odot\mathcal{G}_{2})=
    \Sing(\mathcal{G}_{1})\cap\Sing(\mathcal{G}_{2})\).  In
    particular, if $V\stackrel{\pi}{\longleftarrow} V'$  is permissible for
    \(\mathcal{G}_{1}\odot\mathcal{G}_{2}\), it is also permissible
    for \(\mathcal{G}_{1}\) and for \( \mathcal{G}_{2}\).

    \item Set \(\pi\) as in (1), and let
    \((\mathcal{G}_{1}\odot\mathcal{G}_{2})'\), \(\mathcal{G}'_1\),
    and \(\mathcal{G}'_2\) denote the transforms at \(V'\).  Then:
    \begin{equation*}
    (\mathcal{G}_{1}\odot\mathcal{G}_{2})'=
    \mathcal{G}'_1\odot\mathcal{G}'_2.
    \end{equation*}
\end{enumerate}
}
\end{parrafo}

\section{Idealistic equivalence and integral closure.}

Recall that two ideals, say $I$ and $J$, in a normal domain $R$
have the same integral closure if they are equal for any extension
to a valuation ring (i.e. if $IS=JS$ for every ring homomorphism
$R\to S$ on a valuation ring $S$). The notion extends naturally to
sheaves of ideals. Hironaka considers the following equivalence on
couples $(J,b)$ and $(J',b')$ over a smooth scheme $V$ (see
\cite{Hironaka77}).

\begin{definition}(Hironaka) The couples $(J,b)$
and $(J',b')$ are  {\em idealistic} equivalent on $V$ if $J^{b'}$
and $(J')^b$ have the same integral closure.
\end{definition}

\begin{proposition} Let $(J,b)$ and $(J',b')$ be idealistic
equivalent. Then:

1) $\Sing(J,b)=\Sing(J',b')$.

Note, in particular, that every monoidal transform $V\leftarrow V_1$
on a center $Y\subset \Sing(J,b)=\Sing(J',b')$ defines transforms,
say $(J_1,b)$ and $((J')_1,b')$ on $V_1$.

2)The couples $(J_1,b)$ and $((J')_1,b')$ are  {\em idealistic}
equivalent on $V_1$.
\end{proposition}

If two couples $(J,b)$ and $(J',b')$ are idealistic equivalent over
$V$, the same holds for the restrictions to every open subset of
$V$, and also for restrictions in the sense of \'etale topology,
and even for smooth topology (i.e. pull-backs by smooth morphisms
$W\to V$).

An {\em idealistic exponent}, as defined by Hironaka in
\cite{Hironaka77}, is an equivalence class of couples in the sense
of idealistic equivalence.

\begin{parrafo} {\rm
The previous equivalence relation has an analogous formulation for
Rees algebras, which we discuss below.

}
\end{parrafo}

\begin{definition} Two Rees algebras over $V$, say  $\mathcal{G}=\bigoplus_{k\geq 0}I_kW^k$ and
$\mathcal{G}'=\bigoplus_{k\geq 0}J_kW^k$, are {\em integrally
equivalent}, if both have the same integral closure.

\end{definition}

\begin{proposition} Let $\mathcal{G}$ and $\mathcal{G}'$ be two integrally equivalent Rees algebras over $V$. Then:

1) $\Sing(\mathcal{G})=\Sing(\mathcal{G}')$.

Note, in particular, that every monoidal transform $V\leftarrow V_1$
on a center $Y\subset \Sing(\mathcal{G})=\Sing(\mathcal{G}')$
defines transforms, say $(\mathcal{G})_1$ and $(\mathcal{G}')_1$
on $V_1$.

2)$(\mathcal{G})_1$ and $(\mathcal{G}')_1$ are   integrally
equivalent on $V_1$.
\end{proposition}

If $\mathcal{G}$ and $\mathcal{G}'$ are  {\em integrally}
equivalent on $V$, the same holds for any open restriction, and
also for pull-backs by smooth morphisms $W\to V$.

On the other hand, as $(\mathcal{G})_1$ and $(\mathcal{G}')_1$ are
integrally equivalent, they define the same closed set on $V_1$
(the same singular locus), and the same holds for further monoidal
transformations, pull-backs by smooth schemes, and concatenations
of both kinds of transformations.

\begin{parrafo} \label{inboer}{\rm
For the purpose of resolution problems, the notions of couples and of Rees algebras are equivalent. We first show that any couple can be identified with an algebra, and then show that every Rees algebra arises from a couple. We assign to a couple $(J,b)$ over a smooth scheme $V$ the Rees
algebra, say:$$\mathcal{G}_{(J,b)}=\calo_V[J^bW^b],$$ which is a
graded subalgebra in $\calo_V[W]$. }
\end{parrafo}

\begin{remark}\label{rkboer}
Note that: $\Sing(J,b)=\Sing(\mathcal{G}_{(J,b)})$. In particular, every transformation
\begin{equation*}
\begin{array}{ccccc}
 & V & \stackrel{\pi}{\longleftarrow} & V_{1} & \\ & (J,b) & & (J_1,b) &\\
\end{array}
\end{equation*}
\end{remark}
induces a transformation, say
\begin{equation*}
\begin{array}{ccccc}
 & V & \stackrel{\pi}{\longleftarrow} & V_{1} & \\ & \mathcal{G}_{(J,b)} & & \left( \mathcal{G}_{(J,b)}\right)_1 &\\
\end{array}
\end{equation*} It can be checked that: $\left( \mathcal{G}_{(J,b)}\right)_1=\mathcal{G}_{(J_1,b)}$.

In particular a sequence (\ref{sectransp}) is equivalent to a
sequence (\ref{sectransp11}) over $\mathcal{G}_{(J,b)}$. Moreover, one of them is a resolution if and only if the other is so
(\ref{defresg}).

The following results shows that the class of couples can be embedded
in the class of Rees algebras, in such a way that equivalence
classes are preserved, and that  every Rees algebra is, up to integral equivalence,
of the form $\mathcal{G}_{(J,b)}$ for a suitable $(J,b)$.
\begin{proposition} \label{28pr}Two couples $(J,b)$ and $(J',b')$ are idealistic equivalent
over a smooth scheme $V$, if and only if the Rees algebras
$\mathcal{G}_{(J,b)}$ and $\mathcal{G}_{(J',b')}$ are integrally
equivalent.
\end{proposition}

\begin{proposition} \label{29pr}Every Rees algebra $\mathcal{G}=\bigoplus_{k\geq 0}J_kW^k$, over a smooth scheme $V$,
is integrally
equivalent to one of the form $\mathcal{G}_{(J,b)}$, for a suitable choice of $b$.

\end{proposition}
\proof
Let $U$ be an affine open set in $V$, and assume that the restriction of
$\mathcal{G}$ to $U$ is
$$\mathcal{G}_U=\calo_V(U)[f_1W^{n_1},\dots
,f_sW^{n_s}]= \bigoplus_{k\geq 0}J_k(U)W^k.$$
If $b$ is a common multiple of all positive integers $n_i$, $1\leq i \leq s$, then
$\mathcal{G}_U$ is an finite ring extension of $ \calo_V(U)[J(U)_bW^{b}]$. Finally, since $V$ can be covered
by finitely many affine open sets, we may choose $b$ so that $\mathcal{G}$  is integrally equivalent
to $\mathcal{G}_{(J_b,b)}$.


\section{Diff-algebras, Finite Presentation Theorem, and Koll\'{a}r's tuned ideals.}
\label{sec3}

Here $V$ is smooth over a field $k$, so for each non-negative
integer $s$ there is a locally free sheaf of differential
operators of order $s$, say $Diff^s_k$. There is a natural
identification, say $Diff^0_k=\calo_V$, and for each $s\geq 0$
$Diff^s_k\subset Diff^{s+1}_k$.
We define an
extension of a sheaf of ideals $J\subset \calo_V$, say
$Diff^s_k(J)$, so that over the affine open set $U$,
$Diff^s_k(J)(U)$ is the extension of $J(U)$ defined by adding  $ D(f)$, for all $D\in Diff^s_k(U)$ and $f\in J(U)$.
$Diff^0(J)=J$, and $Diff^s(J)\subset Diff^{s+1}(J)$ as sheaves of
ideals in $\calo_V$. Let $V(I)\subset V$ denote the closed set
defined by an ideal $I\subset \calo_V$.
The order of the ideal $J$ at the local
regular ring $\calo_{V,x}$ is $\geq s$ if and only if $x\in
V(Diff^{s-1}(J))$.

\begin{definition} \label{3def1} We say that a Rees algebra $\bigoplus_{n\geq 0} I_nW^n$, on a smooth scheme $V$,  is a Diff-algebra relative to the field $k$, if: i) $I_n\supset I_{n+1}$ for $n\geq 0$. ii) There is open covering of $V$ by affine open sets $\{U_i\}$,
and for every $D\in Diff^{(r)}(U_i)$, and $h\in I_n(U_i)$, then
 $D(h)\in I_{n-r}(U_i)$ provided $n\geq r$.

Note that (ii) can be reformulated by: ii') $Diff^{(r)}(I_n)\subset I_{n-r}$ for each $n$, and $0\leq r
\leq n$.
\end{definition}
\begin{parrafo}\label{ffbrd}{\em
Fix a closed point $x\in V$, and a regular system of parameters
$\{x_1,\dots,x_n\}$ at $\calo_{V,x}$. The residue field, say $k'$
is a finite extension of $k$, and the completion
$\hat{\calo}_{V,x}=k'[[x_1,\dots , x_n]].$

The Taylor development is the continuous $k'$-linear ring
homomorphism:
$$Tay: k'[[x_1,\dots , x_n]]\to k'[[x_1,\dots , x_n,T_1,\dots ,
T_n]]$$ that map $x_i$ to $x_i+T_i$, $1\leq i \leq n$. So for
$f\in k'[[x_1,\dots , x_n]]$, $Tay(f(x))= \sum_{\alpha \in
\mathbb{N}^n} g_{\alpha}T^{\alpha}$, with $g_{\alpha} \in
k'[[x_1,\dots , x_n]]$. Define, for each $\alpha \in
\mathbb{N}^n$, $\Delta^{\alpha}(f)=g_{\alpha}$. There is a natural
inclusion of $\calo_{V,x}$ in its completion, and it turns out
that
$\Delta^{\alpha}(\calo_{V,x})\subset \calo_{V,x},$
and that $\{ \Delta^{\alpha}, \alpha \in (\nat)^n , 0 \leq
|\alpha| \leq c\}$ generate the $\calo_{V,x}$-module
$Diff^c_k(\calo_{V,x})$ (i.e. generate $Diff^c_k$ locally at $x$).
}
\end{parrafo}
\begin{theorem}\label{thopG} For every Rees algebra $\mathcal{G}$ over a smooth scheme
$V$, there is a Diff-algebra, say $G(\mathcal{G})$ such that:

i) $\mathcal{G}\subset G(\mathcal{G})$.

ii) If $\mathcal{G}\subset \mathcal{G}'$ and $\mathcal{G}'$ is a
Diff-algebra, then $G(\mathcal{G})\subset \mathcal{G}'$.

Furthermore, if $x\in V$ is a closed point, and
$\{x_1,\dots,x_n\}$ is a regular system of parameters at
$\calo_{V,x}$, and if $\mathcal{G}$ is locally generated by
$\mathcal{F}=\{
g_{n_i}W^{n_i},  n_i>0 , 1\leq i\leq m \},$ then
\begin{equation} \small
\mathcal{F'}=\{\Delta^{\alpha}(g_{n_i})W^{n'_i-\alpha}/
g_{n_i}W^{n_i}\in \mathcal{F}, \alpha=(\alpha_1, \alpha_2,\dots ,
\alpha_n) \in (\mathbb{N})^n, \mbox{ and } 0\leq |\alpha| < n'_i
\leq n_i\}
\end{equation}

generates $G(\mathcal{G})$ locally at $x$.
\end{theorem}
(see \cite{VV1}, Theorem 3.4).

\begin{remark}\label{rk33}
1) If $\mathcal{G}_1$ and $\mathcal{G}_2$ are Diff-algebras, then
$\mathcal{G}_1 \odot \mathcal{G}_2$ is also a Diff-algebra.

2)The local description in the Theorem shows that $\Sing(\mathcal{G})=\Sing(G(\mathcal{G}))$.

In fact, as $\mathcal{G}\subset G(\mathcal{G})$, it is clear that
$\Sing(\mathcal{G})\supset \Sing(G(\mathcal{G}))$. For the converse
note that if $\nu_x(g_{n_i})\geq n_i$, then
$\Delta^{\alpha}(g_{n_i})$ has order at least $n_i-|\alpha|$ at
the local ring $\calo_{V,x}$.
\end{remark}

The $G$ operator is compatible with pull-backs by smooth morphisms, and this kind of morphism will arise later
(see \ref{eqdrzo1}). The following Main Lemma, due to Jean Giraud, relates the, say $G$-extensions, with monoidal transformations.

\begin{lemma} (J. Giraud) \label{th41}Let $\mathcal{G}$ be a Rees algebra on a
smooth scheme $V$, and let $V\longleftarrow V_1$ be a permissible
(monoidal) transformation for $\mathcal{G}$. Let $\mathcal{G}_1$
and $G(\mathcal{G})_1$ denote the transforms of $\mathcal{G}$ and
$G(\mathcal{G})$. Then:

1) $\mathcal{G}_1\subset G(\mathcal{G})_1$.

2) $G(\mathcal{G}_1)=G(G(\mathcal{G})_1).$

\end{lemma}
(see \cite{EncVil06}, Theorem 4.1.)



\section{On Hironaka's main invariant.}
\label{sec45}

Hironaka attaches to a couple $(J,b)$ a fundamental invariant for resolution problems, which is a function (see \ref{Agr1}). Here we discuss the role of this function in resolution, and the satellite functions
defined in terms of it. These satellite functions are the main ingredients for the algorithm of resolution in
\cite{Villa89}, for the case of characteristic zero.
\begin{definition}
\label{semicon} Let \(X\) be a topological space, and let \((T, \geq
)\) be a totally ordered set.  A function \(g: X \rightarrow T\) is
said to be {\it upper semi-continuous} if: {\bf i)} \(g\) takes
only finitely many values, and, {\bf ii)} for any $\alpha \in T$
the set $\{x\in X \ / g(x)\geq \alpha \}$ is closed in $X$. The largest value achieved by \(g\) will be denoted by
$\max g$.

We also define
\begin{equation*}
    \vMax g=\{ x\in X : g(x)= \max g
\}
\end{equation*}
 which is a closed subset of \(X\).
\end{definition}

\begin{definition}\label{Agr1} Give a couple $(J,b)$, set
\begin{equation}\label{defeord}
\begin{array}{rccccl}
ord^{}_{(J,b)}: & \Sing(J,b)& \to  \mathbb{Q}\geq 1;&&&ord^{}_{(J,b)}(x)=\frac{\nu_{{J}}(x)}{b}
\end{array}
\end{equation}
where $\nu_{{J}_s}(x)$ denotes the order of
${J}$ at the local regular ring $\mathcal{O}_{V,x}$.

\end{definition}

Note that the function is upper semi-continious; and note also that if $(J_1,b_1)$ and $(J_2,b_2)$ are integrally equivalent, then both functions coincide on $\Sing(J_1,b_1)=\Sing(J_2,b_2)$.

\begin{parrafo}\label{trobo}{\em
Resolution of couples was defined in \ref{parsec} as a composition of permissible transformations, each of which is a monoidal transformation. Every monoidal transformation introduces a smooth hypersurface, and a composition introduces several smooth hypersurfaces. The definition of resolution requires that these hypersurfaces have normal crossings. We define a {\em pair} $(V,E)$ to be  a smooth scheme $V$
together with $E=\{H_1, \dots , H_r\}$ a set of smooth hypersurfaces so that their union has normal crossings. If $Y$ is closed and smooth in $V$ and has normal crossings with $E$ (i.e. with the union of hypersurfaces of $E$), we define a transform of the pair, say  $$(V,E)\leftarrow (V_1,E_1),$$  where $V\leftarrow V_1$ is the blow up at $Y$; and $E_1=\{H_1, \dots , H_r H_{r+1}\}$, where $H_{r+1}$ is the exceptional locus, and each $H_i$ denotes again the strict transform of $H_i$ , for $1\leq i \leq r$.

We define a {\em basic object} to be a pair $(V, E=\{H_1,..,
H_r\})$ together with a couple $(J,b)$ (so
$J_x\neq 0$ at any point $x\in V$). We indicate
this basic object by $$ (V,(J,b),E).$$
If a smooth center $Y$
defines a transformation of  $(V,E)$, and in addition $Y
\subset \Sing(J,b)$, then a transform of the couple $(J,b)$ is
defined. In this case we say that $$ (V,(J,b),E) \longleftarrow
(V_1,(J_1,b),E_1)$$ is a {\em transformation} of the basic object.
A sequence of transformations
\begin{equation}\label{transfuno}
(V,(J,b),E) \longleftarrow (V_1,(J_1,b),E_1)\longleftarrow \cdots
\longleftarrow (V_s,(J_s,b),E_s)
\end{equation}
 is a {\em resolution} of the basic object if
$\Sing(J_s,b)=\emptyset.$

In such case the total transform of $J$ can be expressed as a product, say:
\begin{equation}\label{traboooo}
J\cdot \mathcal{O}_{V_s}=I(H_{r+1})^{c_1}\cdot I(H_{r+2})^{c_2}\cdots
I(H_{r+s})^{c_s
}\cdot J_s
\end{equation}
for some integer $c_i$, where $J_s$ is a sheaf of ideals of order
at most $b-1$, and the  hypersurfaces $H_j$ have normal crossings.

Note that $\{H_{r+1}, \dots , H_{r+s}\}\subset E_s$, and equality holds when $E=\emptyset $. Furthermore, a resolution of a couple $(J,b)$ is attained by a resolution of $(V,(J,b), E=\emptyset)$ (see \ref{parsec}).

}
\end{parrafo}

\begin{parrafo} {\bf The first satellite functions.}\label{lasat1} {\em (see 4.11, \cite{EncVil97:Tirol}) Consider, as above,  transformations
\begin{equation}\label{Atransfuno}
(V,(J,b),E) \longleftarrow (V_1,(J_1,b),E_1)\longleftarrow \cdots
\longleftarrow (V_s,(J_s,b),E_s)
\end{equation}
which is not necessarily a resolution, and let $\{H_{r+1}, \dots , H_{r+s}\}(\subset E_s)$ denote the exceptional hypersurfaces introduced by the sequence of blow-ups. We may assume, for simplicity that these hypersurfaces are irreducible. There is a well defined factorization of the sheaf of ideals $J_s \subset \calo_{V_s}$, say:
\begin{equation}\label{eqfcword}
J_s=I(H_{r+1})^{b_1}I(H_{r+2})^{b_2}\cdots I(H_{r+s})^{b_s}\cdot
\overline{J}_s
\end{equation} so that $\overline{J}_s$ does not vanish along
$H_{r+i}$, $1\leq i \leq s$.

Define  $\vword^{d}_{(J_s,b)}$ (or simply $\vword^{d}_{s}$):
\begin{equation}\label{eqdword}
\begin{array}{rccccl}
\vword^{d}_s: & \Sing(J_s,b)& \to  \mathbb{Q};&&&\vword^d_s(x)=\frac{\nu_{\overline{J}_s}(x)}{b}
\end{array}
\end{equation}
where $\nu_{\overline{J}_s}(x)$ denotes the order of
$\overline{J}_s$ at $\calo_{V_s,x}$.  It has the following properties:

 1) The function is
upper semi-continuous. In particular $\vMax  \vword$ is closed.

2) For any index $ i \leq s$, there is an expression
$$J_i=I(H_1)^{b_1}I(H_2)^{b_2}\cdots I(H_i)^{b_i}\cdot
\overline{J}_i,$$ and hence a function $\vword_i:
\Sing(J_i,b) \to {\mathbb Q}$ can also be defined.

3)  If each transformation of basic objects  $(V_i,(J_i,b_i), E_i) \leftarrow (V_{i+1},(J_{i+1},b_{i+1} )E_{i+1})$ in
(\ref{Atransfuno}) is defined with center $Y_i \subset \vMax \vword_i $,
then $$ \max \vword \geq \max \vword_1\geq \dots \geq \max
\vword_s.$$

}
\end{parrafo}

\begin{parrafo}\label{rk55} {\em  Let us stress here on the fact that the previous definition of the function  $\vword $ (\ref{eqdword}) (and the factorization in
(\ref{eqfcword})), grow from the function introduced in Def \ref{Agr1}.
Fix $H_{r+i}$ as in (\ref{eqfcword}), and define a function $exp_i$ along the points in $\Sing(J_s,b)$  by setting $exp_i(x)=\frac{b_i}{b}$ if $x\in H_{r+i}\cap \Sing(J_s,b)$, and $exp_i(x)=0$ otherwise.

Induction on the integer $s$ allow us to express each rational number $exp_i(x)$ in terms of the functions
$ord_{(J_{s'},b)}$, for $s'<s$ (in terms of these functions $ord_{(J_{s'},b)}$  evaluated at the generic points, say $y_{s'}$, of the centers
 $Y_{s'} ( \subset V_{s'})$ of the monoidal transformation). Finally note that
$$\vword^d_{(J_s,b)}(x)=ord_{(J_s,b)}(x)- exp_1(x)-exp_2(x)-\dots - exp_s(x).$$
Therefore  all these functions grow from Hironaka's functions $ord_{(J_i,b)} $, so we call them "satellite functions" (p.p. 187 \cite{EncVil97:Tirol}). In particular, if $(J,b)$ and $(J',b')$ are idealistic equivalent at $V$, then (\ref{transfuno}) induces  transformations
of $(V,(J',b'),E)$; moreover $\Sing(J_s,b)=\Sing(J'_s,b')$, and $\vword_{(J_s,b)}=\vword_{(J'_s,b')}$ as functions (and the exponent functions $exp_i$ coincide).

The general strategy to obtain resolution of a basic objects
$(V, (J,b), E) $ (and hence of couples $(J,b)$), is to produce a sequence of transformations as in
(\ref{Atransfuno}), so that  $\overline{J}_s=\calo_{V_s}$ in an open neighborhood of $\Sing(J_s,b)$
(\ref{eqfcword}). This amounts to saying that  $\vword^d_s(x)=0$ at any $x\in \Sing( J_s,b)$, or equivalently, that $\max \vword=0$.
When this condition holds we say that the transform  $(V_s,(J_s,b_s), E_s) $ is in the {\em monomial case}.

If this condition is achieved, then we may assume
$J_s=I(H_{r+1})^{b_1}I(H_{r+2})^{b_2}\cdots I(H_{r+s})^{b_s}$, and it is simple to extend, in this case, sequence (\ref{Atransfuno}) to a resolution. In fact one can extend it to a resolution by choosing centers as intersections of the exceptional hypersurfaces; which is a simple combinatorial strategy defined in terms of the exponents $b_i$.

The point is that there is a particular kind of basic object, which we define below, called
{\em simple basic objects}, with the following property:  if we know how to produce resolution of simple basic objects then we can define (\ref{Atransfuno}) so as to achieve the monomial case.
The point is that resolution of simple basic objects should be achieved by some form of induction. This is why we call them simple.

\begin{definition} Let $V$ be a smooth scheme, and let $(J,b)$ be a couple.

1) $(J,b)$ is said to be a {\em simple couple} if $ord^{}_{(J,b)}(x)=1$ for any $x \in \Sing(J,b)$ (see
 \ref{defeord}).

 2) $(V,(J,b), E)$ is a {\em simple basic object} if $(J,b)$ is simple and $E=\emptyset$.

\end{definition}}
\end{parrafo}

\begin{parrafo} {\bf  Second satellite function: the inductive function t.} \label{poift}{\em
 (See 4.15 ,\cite{EncVil97:Tirol}.)
Consider
\begin{equation}
\label{esol6}
\begin{array}{cccccccc}
(V,(J,b),E)& \leftarrow &(V_1,(J_1,b),E_1) & \leftarrow &\cdots
&\leftarrow & (V_s,(J_s,b),E_s),&
\\
\end{array}
\end{equation}
as before, where each $V_i\leftarrow V_{i+1}$ is defined with center $Y_i
\subset \vMax \vword_i$, so that:
\begin{equation}\label{eqdgw}
 \max \vword \geq \max
\vword_1\geq \dots \geq \max \vword_s.
\end{equation}

We now define a function $t_s$, only under the assumption that $\max
\vword_s>0$. In fact, $\max \vword_s=0$ in the monomial case, which is
easy to resolve.

Set $s_0 \leq s$ such that
$$ \max \vword \geq \dots \geq \max \vword_{s_{0-1}}> \max
\vword_{s_0}= \max \vword_{s_{0+1}}= \dots= \max \vword_{s},$$ and set:
$$ E_s=E_s^+ \sqcup E_s^- \mbox{ (disjoint union)},$$ where
$E_s^-$ are the strict transform of hypersurfaces in $E_{s_0}$.
Define
$$ t_s: \Sing(J_s,b) \to {\mathbb Q} \times {\mathbb N}  \mbox{ (ordered lexicographically). }$$
$$ t_s(x)=(\vword_s(x), n_s(x))\hskip 1cm  n_s(x)= \sharp \{H_i\in E_s^-  | x\in H_i\}$$
where $\sharp S$ denotes the total number of elements of a set $S$. One can check that:

i) the function is upper semi-continuous. In particular $\vMax t_s$
is closed.

ii) There is a function $t_i$ for any index $ i \leq s$; and if $(J,b)$ and $(J',b')$ are integrally equivalent over $V$, then (\ref{esol6}) induces a sequence of transformations of $(V,(J',b'),E)$ and the corresponding functions $t_i$ coincide along $\Sing(J_i,b)=\Sing(J'_i,b')$.

iii) If each $(V_i,E_i) \leftarrow (V_{i+1},E_{i+1})$ in
\ref{esol6} is defined with center $Y_i \subset \vMax t_i $, then
\begin{equation}\label{eqdgt}
 \max t \geq \max t_1\geq \dots \geq \max t_s.
\end{equation}
iv) If $\max t_s=( \frac{d}{b}, r)$ (here $\max
\vword_s=\frac{d}{b}$) then $\vMax t_s \subset \vMax \vword_s$.

Recall that the functions $t_i$ are defined only if $\max \vword_i>0$. We say that a sequence of transformations is $t$-{\em permissible} when condition iii) holds; namely when $Y_i \subset \vMax t_i $ for all $i$.
}
\end{parrafo}
\begin{definition}\label{defbirth}  Consider, as above, a sequence
\begin{equation}
\label{reesol6}
\begin{array}{cccccccc}
(V,(J,b),E)& \leftarrow &(V_1,(J_1,b),E_1) & \leftarrow &\cdots
&\leftarrow & (V_s,(J_s,b),E_s),&
\\
\end{array}
\end{equation}
so that $Y_i \subset \vMax \vword_i $
for $0\leq i\leq s$, and furthermore: that $Y_i \subset \vMax t_i (\subset \vMax \vword_i)$ if $\max  \vword_i>0$. The decreasing sequence of values  (\ref{eqdgw}) will hold, and that also (\ref{eqdgt}) holds if
$\max  \vword_s>0$. We now attach an index, say $r$, to the basic object $(V_s,(J_s,b),E_s)$, defined in terms of the sequence of transformations.

i) If $\max  \vword_s>0$ set $r$ as the smallest index so that
$\max  t_r=\max  t_{r+1}=\cdots = \max  t_s.$

ii) If $\max  \vword_s=0$ set $r$ as the smallest index so that
$\max  \vword_r=0$.

\end{definition}
\begin{parrafo}\label{comnt} {\em The satellite functions were defined for suitable sequences of transformations of basic objects. The main properties of the Inductive Function $t$ can be stated as follows:

1) There is a simple basic object naturally attached to the function.

2) This simple basic object can be chosen so as to be well defined up idealistic equivalence.

The following Proposition clarifies  1), whereas 2) will be addressed later (see \ref{parKol}).
}
\end{parrafo}
\begin{proposition} \label{propdegbo} Assume that a sequence of $s$ transformations  (\ref{reesol6})  of $(V,(J,b),E)$  is defined in the same conditions as above, and that $\max  \vword_s>0$.  Fix $r$ as in Def  \ref{defbirth} , i).

There is a simple couple $(J'_r,b')$ at $V_r$, so that the simple basic object $(V_r,(J'_r,b'),E'_r=\emptyset)$ has the following property:

Any sequence of transformations of $(V_r,(J'_r,b'),E'_r=\emptyset)$, say
\begin{equation}
\label{resolll6}
\begin{array}{cccccccc}
(V_r,(J'_r,b'),\emptyset)& \leftarrow &(V'_{r+1},(J'_{r+1},b'),E'_{r+1}) & \leftarrow &\cdots
&\leftarrow & (V'_S,(J'_S,b'),E'_S),&
\\
\end{array}
\end{equation}
induces a $t$-permissible sequence of transformation of the basic object $(V_r,(J_r,b),E_r)$, say:
\begin{equation}
\label{resollll6}
\begin{array}{cccccccc}
(V_r,(J_r,b),E_r)& \leftarrow &(V'_{r+1},(J_{r+1},b),E_{r+1}) & \leftarrow &\cdots
&\leftarrow & (V'_S,(J_S,b),E_S),&
\\
\end{array}
\end{equation}
by blowing up at the same centers, with the following condition on the functions $t_j$ defined for this last sequence (\ref{resollll6}):

a) $\vMax t_k=\Sing(J'_k,b')$, $r\leq k \leq S-1$.

b) $\max  t_r=\max  t_{r+1}=\cdots = \max  t_{S-1} \geq \max  t_{S}.$

c) $\max  t_{S-1} = \max  t_{S}$ if and only if $\Sing(J'_S,b')\neq \emptyset$, in which case $\vMax t_S=\Sing(J'_S,b')$.

\end{proposition}
\proof (see \ref{ammr})

\begin{remark} \label{rrr43}So if we take (\ref{resolll6})  to be a resolution of  $(V_r,(J'_r,b'),\emptyset)$, we can extend the first $r$ steps of (\ref{reesol6}), say
\begin{equation}
\label{ressol6}
\begin{array}{cccccccc}
(V,(J,b),E)& \leftarrow &(V_1,(J_1,b),E_1) & \leftarrow &\cdots
&\leftarrow & (V_r,(J_r,b),E_r),&
\\
\end{array}
\end{equation} with the transformations of sequence (\ref{resollll6}), say:
\begin{equation*}
\begin{array}{ccccccccc}
(V,(J,b),E)& \leftarrow&\cdots &(V_r,(J_r,b),E_r) & \leftarrow(V'_{r+1},(J_{r+1},b),E_{r+1})  &\cdots
&\leftarrow & (V'_S,(J_S,b),E_S),&
\\
\end{array}
\end{equation*}
and now
$\max  t_r=\max  t_{r+1}=\cdots = \max  t_{S-1} > \max  t_{S}.$

In other words, the Proposition asserts that if we know how to define resolution of simple basic objects, then we can force the value $\max t$ to drop. Note that for a fixed basic object $(V,(J,b),E)$ there are only finitely many possible values of $\max t$ in any sequence of permissible monoidal transformations. As indicated in \ref{rk55}, resolution of simple basic objects would lead to the case $\max \vword_S=0$, also called the  monomial case, which is easy to resolve.
The conditions of Prop. \ref{propdegbo} hold for $s=0$, namely when (\ref{reesol6})  is simply $(V,(J,b),E)$. So given  $(V,(J,b),E)$, this Proposition indicates how to define a sequence of transformations that takes it to the monomial case (provided we know how to resolve simple basic objects). Moreover, a unique procedure of resolution of simple basic objects would define, for each  $(V,(J,b),E)$, a unique sequence of transformations that takes it to the monomial case.

\end{remark}

\begin{remark} \label{ammr}
 A general property of simple couples is that any transform, say $(J_1,b)$, is again simple. An outstanding  property of the satellite functions is that they are upper semi-continuous, and a simple basic object can be attached to the highest value achieved by the function. Let the setting and notation be as in \ref{lasat1}. If $\max \vword_i=\frac{d}{b}$, and $d\geq b$, set
\begin{equation}\label{ecdwrd1}
(J_i'',b'')= (\overline{J}_i,d)
\end{equation}
If $\max \vword_i=\frac{d}{b}$, and $d< b$, set
\begin{equation}\label{ecdwrd2}
(J_i'',b'')= (J^{d}+ (\overline{J}_i)^b,bd).
\end{equation}
One can check that if $\max \vword_i=\frac{d}{b}$, then $\vMax
\vword_i=\Sing(J_i'',b'')$. So the points where the functions takes its highest values is the
singular locus of a simple basic object. Furthermore, this link is preserved by transformations in the following sense. Assume, as in \ref{poift} that a sequence (\ref{esol6})
is defined by centers $Y_i\subset \vMax \vword_i$, and set $s_0$ as the smallest index so that $\max \vword_{s_0}=\max \vword_{s}$. One can check that, for each index $i\geq s_0$: $(J_{i+1}'',b'')$ is the transform of $(J_{i}'',b'')$.

A similar property holds for for the inductive function $t$. In fact Proposition \ref{propdegbo} establishes
an even stronger link of the value $\max t$ with a simple basic.
Given a positive integer $h$, let $G_h$ be the set
of all subsets $\mathcal{F}\subset E^-_i$,
$\mathcal{F}=\{H_{j_1},\dots ,H_{j_h}\}$ (with $h$ hypersurfaces).
For each positive integer $m$ define
$\mathcal{H}_h(m)=\prod_{\mathcal{F}\in G_h} \sum_{{H_{i_j}}\in
\mathcal{F}} I(H_{i_j})^m.$

Set $\max t_i=(\frac{d}{b},h)$. If $d\geq b$ set $$(J'_i,b')=(J''_i+\mathcal{H}_h(d), d),$$ with
$J''_i$ as in (\ref{ecdwrd1}). If $d< b$ set $$(J_i',b')=(J''_i+\mathcal{H}_h(bd), bd),$$ with
$J''_i$ as in (\ref{ecdwrd2}).
 Note also the $(J_i',b')$ is simple, and
$\Sing(J_i',b')=\vMax t_i \subset \Sing(J_i',b')=\vMax \vword_i$.

One can check that $(J'_r, b')$ fulfills the condition in Prop \ref{propdegbo} (see Th 7.10, \cite{EncVil97:Tirol}).
\end{remark}

\begin{parrafo} \label{elpdrob}{\bf New operations on basic objects.}  {\em There are two natural operations on basic objects, which we discuss below, that play a central role in Hironaka's definition of {\em  invariance} of the main function introduced in Definition  \ref{Agr1}, and also for the proof of the second property stated in \ref{comnt}.  Recall the definition of a basic object over a smooth scheme $V$  over a field $k$ in \ref{trobo}, say
$ (V,(J,b),E)$, where a {\em pair} $(V,E)$ is defined by $E=\{H_1, \dots , H_r\}$, a set of smooth hypersurfaces with normal crossings, and  $(J,b)$ is a couple on $V$. Let now
\begin{equation}\label{eqdrtbo}
V \stackrel{\pi}{\longleftarrow} U
\end{equation}
be defined by:

A) An open set $U$ of $V$ (Zariski or \'etale  topology).

B) The projection of $U=V\times \mathbb{A}_k^n$ on the first coordinate, where $ \mathbb{A}_k^n$ denotes the $n$-dimensional affine scheme, for some positive integer $n$.

In both cases there is a natural notion of pull-backs of the pair $(V,E)$  to $(U, E_U)$, where $E_U$ consists of the pull-back of hypersurfaces in $E$. There is also a notion of pull-back of the couple $(J,b)$, say $(J_U,b)$, by restriction in case A), and taking the total transform in case B).

In this way  we attach to (\ref{eqdrtbo}) a  {\em pull-back} of basic objects, say
\begin{equation}\label{eqdrtbo1}
 (V,(J,b),E) \stackrel{\pi}{\longleftarrow}  (U,(J_U,b),E_U)
\end{equation}

The first observation is that the singular locus and Hironaka«s function in Definition \ref{Agr1} are compatible with pull-backs; i.e.
$\Sing(J_U,b)=\pi^{-1}(\Sing(J,b))$, and for $x\in \Sing(J_U,b)$:
\begin{equation}\label{eqdjtbo6}
 \ord^{}_{(J_U,b)}(x)=\ord^{}_{(J,b)}(\pi(x)).
\end{equation}
A similar compatibility property holds for satellite functions. Recall that these functions
were defined for transformations of a basic object, say
\begin{equation}\label{ABtranso}
(V,(J,b),E) \longleftarrow (V_1,(J_1,b),E_1)\longleftarrow \cdots
\longleftarrow (V_s,(J_s,b),E_s)
\end{equation}
to which we attached a well defined factorization of the sheaf of ideals $J_s \subset \calo_{V_s}$, say:
\begin{equation}\label{eqword1}
J_s=I(H_{r+1})^{b_1}I(H_{r+2})^{b_2}\cdots I(H_{r+s})^{b_s}\cdot
\overline{J}_s
\end{equation} so that $\overline{J}_s$ does not vanish along
$H_{r+i}$, $1\leq i \leq s$; where  $\{H_{r+1}, \dots , H_{r+s}\}(\subset E_s)$ denote the exceptional hypersurfaces introduced by the sequence of blow-ups. Let now
$$ (V_s,(J_s,b),E_s) \stackrel{\pi}{\longleftarrow}  (U_s,((J_s)_U,b),(E_s)_U)$$
be a pull-back of $ (V_s,(J_s,b),E_s) $.  Note that $\{\pi^{-1}(H_{r+1}), \dots , \pi^{-1}(H_{r+s})\}(\subset (E_s)_U)$, indicate the non-smooth locus of the composite morphism $V\leftarrow (U_s)$, and there is a
natural pull-back of (\ref{eqword1}), say
\begin{equation}\label{eqword2}
(J_s)_U=(I(H_{r+1}))_U^{b_1}(I(H_{r+2}))_U^{b_2}\cdots (I(H_{r+s}))_U^{b_s}\cdot
(\overline{J}_s)_U.
\end{equation}
}
\end{parrafo}

This shows that the function $\vword$ in (\ref{eqdword}) is also compatible with pull-backs.
One can go one step further:

\begin{definition} \label{lsotr} Define a {\em local sequence of transformations} of a basic object
$(V,(J,b), E)$, to be a sequence
\begin{equation}\label{ABtranso}
(V,(J,b),E) \longleftarrow (V'_1,(J_1,b),E_1)\longleftarrow \cdots
\longleftarrow (V'_s,(J_s,b),E_s),
\end{equation}
where $(V,(J,b),E) \longleftarrow (V'_1,(J_1,b),E_1)$, and each
$(V'_i,(J_i,b),E_i) \longleftarrow (V'_{i+1},(J_{i+1},b),E_{i+1})$,
 is a pull-back, or a pull-back followed by a usual transform of a basic object (as  defined in \ref{trobo}).
\end{definition}

\begin{parrafo}\label{exacdlt} {\em
The same discussion in \ref{lasat1} applies for a local sequence of transformations: if each
\begin{equation}\label{eqaff}
(V'_i,(J_i,b),E_i) \longleftarrow (V'_{i+1},(J_{i+1},b),E_{i+1})
\end{equation}
is a pull-back followed by the blow up at a center $Y_i\subset \vMax \vword_i$ , then
$$ \max \vword \geq \max \vword_1\geq \dots \geq \max
\vword_s.$$
A similar argument applies  for the function $t$ : if the previous condition holds, and $\max \vword_s>0$, then the functions $t_i$ can be defined for $0\leq i \leq s$. Furthermore, if each "local transformation"  (\ref{eqaff}) is the blow up at a center $Y_i\subset \vMax t_i$ followed by a pull-back, then
$$ \max t \geq \max t_1\geq \dots \geq \max
t_s.$$
In this case we shall say that (\ref{ABtranso}) is a $t$-permissible sequence of  ''local transformations".
}
\end{parrafo}
\begin{parrafo}\label{cftr}{\em The following Proposition is stronger then Proposition \ref{propdegbo}, in fact it expresses a stronger property of the inductive function $t$. As we shall see later, after discussing a weaker equivalence notion in the next section, this stronger version will ensure statement 2) in \ref{comnt}.
}
\end{parrafo}

\begin{proposition} \label{prdegbo1} Let the setting and notation be as above, where
   (\ref{ABtranso}) is  a $t$-permissible sequence of local transformations of $(V,(J,b),E)$, and $\max  \vword_s>0$. Let $r$ be the smallest index so that  $\max t_r=\max t_s$. There is a simple basic object $(V_r,(J'_r,b'),E'_r=\emptyset)$ with the following property.
  An arbitrary local sequence of transformations of $(V_r,(J'_r,b'),E'_r=\emptyset)$, say
\begin{equation}
\label{reesolll6}
\begin{array}{cccccccc}
(V_r,(J'_r,b'),\emptyset)& \leftarrow &(V'_{r+1},(J'_{r+1},b'),E'_{r+1}) & \leftarrow &\cdots
&\leftarrow & (V'_S,(J'_S,b'),E'_S),&
\\
\end{array}
\end{equation}
induces a $t$-permissible local sequence of transformation of the basic object $(V_r,(J_r,b),E_r)$, say:
\begin{equation}
\label{reesollll6}
\begin{array}{cccccccc}
(V_r,(J_r,b),E_r)& \leftarrow &(V'_{r+1},(J_{r+1},b),E_{r+1}) & \leftarrow &\cdots
&\leftarrow & (V'_S,(J_S,b),E_S),&
\\
\end{array}
\end{equation}
with the following condition on the functions $t_j$ defined for this last sequence (\ref{reesollll6}):

a) $\vMax t_k=\Sing(J'_k,b')$, $r\leq k \leq S-1$.

b) $\max  t_r=\max  t_{r+1}=\cdots = \max  t_{S-1} \geq \max  t_{S}.$

c) $\max  t_{S-1} = \max  t_{S}$ if and only if $\Sing(J'_S,b')\neq \emptyset $, in which case $\vMax t_S=\Sing(J'_S,b')$.
\proof The couple $(J'_r, b')$, defined in \ref{ammr} , also fulfills these conditions.

\end{proposition}

\section{A weaker equivalence notion.}
\label{sec46}
The real strength of Hironaka's function in Definition \ref{Agr1}, and hence of the satellite functions  in
\ref{lasat1} and \ref{poift}, cannot be understood unless we discuss a weak form of equivalence on couples $(J,b)$, which we do in this section (see Definition 6.15, \cite{EncVil97:Tirol}).

\begin{parrafo}\label{equifte}
{\em Fix a smooth scheme $V$ and a couple $(J,b)$.  Note that the closed set attached to it coincides with that attached to $(J^2,2b)$; namely $\Sing(J,b)=\Sing(J^2,2b)$.
The same holds if we take a pull-back (\ref{eqdrtbo}) either of type A) or B); and the same holds after any local sequence of transformations of $(V,(J,b),\emptyset)$ and $(V,(J^2,2b),\emptyset)$.

A similar property holds for two couples on $V$ which are idealistic equivalent, say $(J_1,b_1)$ and
$(J_2,b_2)$. In fact any pull-back defines two idealistic equivalent couples, and any monoidal transformation of idealistic couples  remain idealistic equivalent.

\begin{definition}\label{deeqft}
Two couples $(J^{(i)},b_i)$ , $i=1,2$, or two basic objects $(V,J^{(i)},b_i);,E)$ , $i=1,2$, are said to be {\em weakly equivalent} if:
$\Sing(J^{(1)},b_1)=\Sing(J^{(2)},b_2),$ and if any local sequence of transformations of one of them:
\begin{equation*}
(V',(J^{(i)},b_i),E') \longleftarrow (V'_1,(J^{(i)}_1,b),E'_1)\longleftarrow \cdots
\longleftarrow (V'_s,(J^{(i)}_s,b),E_s),
\end{equation*}
defines a local sequence of transformations of the other, and $\Sing(J_i^{(1)},b_1)=\Sing(J_i^{(2)},b_2)$, $0\leq i\leq s$.
\end{definition}
}
\end{parrafo}
\begin{parrafo}\label{equift11}
{\em
It is important to point out that a first example of weak equivalence appeared already in Proposition \ref{prdegbo1}: if
$(V_r,(J'_r,b'),\emptyset)$ and $(V_r,(J''_r,b''),\emptyset)$ are two basic objects which fulfill the property of that Proposition, then they are weakly equivalent.

Note that if $(V,(J^{(i)},b_i), E)$ , $i=1,2$, are weakly equivalent as above, then also their transforms
$(V_s,J_s^{(i)},b_i);,E_s)$   are weakly equivalent, for $i=1,2$. So this equivalence is preserved after an arbitrary  local sequence of transformations.

Note also that integral equivalence implies weak equivalence. Hironaka proves a suitable convers which we will address later. This converse will clarify the second property in \ref{comnt}.

}\end{parrafo}

\begin{lemma} {\bf Hironaka's First Main Lemma.} \label{HFML1}{\em (See \cite{Hironaka77};  or 7.1, \cite{EncVil97:Tirol}.)}  If two basic objects $(V,(J^{(i)},b_i);,E)$, $i=1,2$, are weakly equivalent, then
$$ \ord_{(J^{(1)},b_1)}(x) =\ord_{(J^{(2)},b_2)}(x)$$
for all $x\in \Sing(J^{(1)},b_1)=\Sing(J^{(2)},b_2)$.
\end{lemma}

\begin{parrafo}\label{equift12}
{\em  Assume that $(V,(J^{(i)},b_i), E)$ , $i=1,2$, are weakly equivalent, and consider a sequence of monoidal transformations
\begin{equation*}
(V,(J^{(1)},b),E) \longleftarrow (V_1,(J^{(1)}_1,b),E_1)\longleftarrow \cdots
\longleftarrow (V_s,(J^{(1)}_s,b),E_s).
\end{equation*}
In \ref{lasat1} we attached to such data an expression
$J^{(1)}_s=I(H_{r+1})^{c_1}I(H_{r+2})^{c_2}\cdots I(H_{r+s})^{c_s}\cdot
\overline{J}^{(1)}_s$
so that $\overline{J}^{(1)}_s$ does not vanish along
$H_{r+i}$, $1\leq i \leq s$.
The previous Proposition, and the discussion in \ref{rk55}, assert that the same monoidal transformations define
\begin{equation}
(V,(J^{(2)},b),E) \longleftarrow (V_1,(J^{(2)}_1,b),E_1)\longleftarrow \cdots
\longleftarrow (V_s,(J^{(2)}_s,b),E_s),
\end{equation}
 an expression
$J^{(2)}_s=I(H_{r+1})^{d_1}I(H_{r+2})^{d_2}\cdots I(H_{r+s})^{d_s}\cdot
\overline{J}^{(2)}_s,$
and for any $x\in \Sing(J_s^{(1)},b_1)=\Sing(J_s^{(2)},b_2)$: $$\vword^{(1)}(x)=\vword^{(2)}(x),\mbox{ and }
exp_1^{(1)}(x)=exp_2^{(2)}(x)$$ for $i=1,\dots, s$ (see  \ref{rk55}).
Similar equalities hold for the function $t$ in  \ref{poift}, and for $t$-permissible transformations.

Moreover, the discussion in \ref{exacdlt} show that these equalities of satellite functions also extends to the case of a local sequence of transformation.

}\end{parrafo}
\begin{parrafo}\label{trnvr1}{\em  Let $V$ be a smooth scheme, so the local ring at a closed point  $\calo_{V,x}$ is regular. Then the associated graded ring, say $gr_M(\calo_{V,x})$ is a polynomial ring, and
$Spec(gr_M(\calo_{V,x}))=\mathbb{T}_{V,x}$ is the tangent space at $x$, which is a vector space.
Fix a vector space $\mathbb{W}$. A vector $v\in \mathbb{W}$ defines a translation, namely
$ tr_v(w)=v+w$ for $v\in \mathbb{W}$.
An homogeneous ideal in the polynomial ring $gr_M(\calo_{V,x})$ defines a closed set , say $C$, in the vector space $\mathbb{T}_{V,x}$ , which is a union of lines through the origin. There is a largest linear subspace, say $L_C$, so that $ tr_v(C)=C$ for any $v\in L_C$.
If we take, for example, $C$ to be defined by $XY$ in the polynomial ring $k[X,Y,Z]$, then
$L_C$ is the $Z$-axis in $\mathbb{A}^3$. We discuss in \ref{slsdh} how equations defining the subspace $L_C$ arise.

Let  $(J,b)$ be a couple on $V$, and fix $x\in \Sing(J,b)$. Hironaka considers the closed set, say  $C_{(J,b)}$ in $\mathbb{T}_{V,x}$  defined by the ideal spanned by $In_b(J) (\subset M^b/M^{b+1})$; and then he defines the integer  $\tau_{(J,b)}(x)$ to be the codimension of the subspace $L_C$ (in $\mathbb{T}_{V,x}$). An important property of this subspace is that for any smooth center $Y$ in $V$, containing the point $x$ and included in $\Sing(J,b)$, the tangent space, say  $\mathbb{T}_{Y,x}$, is a subspace of $L_C$ (as  subspace of  $\mathbb{T}_{V,x}$). Note that  $ord_{(J,b)}(x)>1$ iff $L_C=\mathbb{T}_{V,x}$ (iff $\tau=0$).

}\end{parrafo}

\begin{lemma} {\bf Hironaka's Second Main Lemma.} \label{HFML2}{\em (See \cite{AHV1})} If two basic objects $(V,(J^{(i)},b_i);,E)$ , $i=1,2$, are weakly equivalent, then
$$ \tau_{(J^{(1)},b_1)}(x) =\tau_{(J^{(2)},b_2)}(x)$$
for $x\in \Sing(J^{(1)},b_1)=\Sing(J^{(2)},b_2)$.

\end{lemma}

\begin{parrafo}\label{agre1}{\bf From Couples to Rees algebras.}
{\em Given a couple $(J,b)$ on a smooth scheme $V$, a function  $ord_{(J,b)}:\Sing(J,b)\to \mathbb{Q}$ was defined in Def. \ref{Agr1} .
In \ref{inboer} we show that every couple defines a Rees algebra $\mathcal{G}_{(J,b)}$, and this assignment is such that $\Sing(J,b)=\Sing(\mathcal{G}_{(J,b)})$. Moreover, the assignment is compatible with transformations (see \ref{rkboer}).

We reformulate Hironaka's  function on the class of Rees algebras, which is the analog
to that defined for couples. Fix $\mathcal{G}=\bigoplus_{k\geq 1}I_kW^k$ and $x\in
\Sing(\mathcal{G})$. Given $f_nW^n\in I_nW^n$ set
$ord_{x,n} (f_n)=\frac{\nu_x(f_n)}{n}\in \mathbb{Q};$ called the
order of $f_n$ (weighted by $n$), where $\nu_x$ denotes the order
at the local regular ring $\calo_{V,x}$. As $x\in
\Sing(\mathcal{G})$ it follows that $ord_{x,n} (f_n)\geq 1.$ Define $$ord_{\mathcal{G}}(x)= inf_{}\{ord_{x,n} (f_n); f_nW^n\in
I_nW^n\}.$$

So $ord_{\mathcal{G}}(x)= inf_{n\geq 1}\{\frac{\nu_x(I_n)}{n}\}$, and $ord_{\mathcal{G}}(x)\geq 1$ for every $x\in
\Sing(\mathcal{G})$. }
\end{parrafo}
\begin{proposition}\label{3propsingZ}
1)  Let $\mathcal{G}$  be a Rees algebra generated  by
$\mathcal{F}=\{ g_{n_i}W^{n_i}, n_i>0 , 1\leq i\leq m \}$  over $V$, then
$ord_{\mathcal{G}}(x)= inf_{}\{ord_{x,n_i} (g_{n_i}); 1\leq i\leq m\}.$

2) If $\mathcal{G}$ and $ \mathcal{G}'$ are Rees algebras with
the same integral closure (e.g. if $\mathcal{G}\subset
\mathcal{G}'$ is a finite extension), then, for every $x\in
\Sing(\mathcal{G})(=\Sing(\mathcal{G}'))$
$$ord_{\mathcal{G}}(x)=ord_{\mathcal{G}'}(x).$$

3) If $\mathcal{G}=\mathcal{G}_{(J,b)}$,  then $ord_{\mathcal{G}}=ord_{(J,b)}$ as functions on
$\Sing(J,b)=\Sing(\mathcal{G}_{(J,b)})$. (So if $N$ is a common multiple of all $n_i , 1\leq i\leq m$,
then $ord_{\mathcal{G}}(x)=inf_{}\{\frac{\nu_x(I_N)}{N}\}$).
\end{proposition}


\begin{parrafo}\label{ddfstll}{\em
We have defined satellite function on basic objects by considering, for a sequence
\begin{equation}
\label{resol341pp}
\begin{array}{cccccccc}
(V,(J,b),E)& \leftarrow &(V_1,(J_1,b),E_1) & \leftarrow & \dots
&\leftarrow & (V_s,(J_s,b),E_s),&
\\
\end{array}
\end{equation}
a natural  factorization:
\begin{equation}\label{resol3-1}
J_s=I(H_1)^{b_1}I(H_2)^{b_2}\cdots I(H_s)^{b_s}\cdot
\overline{J}_s
\end{equation} so that $\overline{J}_s$ does not
vanish along $H_i$, $1\leq i \leq s$. Recall that $\vword_s$ is defined in terms of this expression
(see  (\ref{eqdword})).
Every Rees algebra $\mathcal{G}$ can be identified, up to integral closure, with one of the form $\mathcal{G}=\mathcal{G}_{(J,b)}$.
Consider now the Rees algebra
$\mathcal{G}=\mathcal{G}_{(J,b)}$, then (\ref{resol341pp}) induces:
\begin{equation}\label{sectransp2}
\begin{array}{cccccccc}
 & (V,E) & \stackrel{\pi_1}{\longleftarrow} &( V_{1},E_1) & \stackrel{\pi_2}{\longleftarrow}
 &\ldots& \stackrel{\pi_k}{\longleftarrow}&(V_s, E_s).\\ & \mathcal{G} & & (\mathcal{G})_1 &&&&(\mathcal{G})_s\\
\end{array}
\end{equation}
In this case  $(\mathcal{G})_i=\mathcal{G}_{(J_i,b)}$ (\ref{rkboer}), so $\Sing((J_i,b) )=\Sing(\mathcal{G}_i )$, and  function
  \begin{equation}
 \vword_{\mathcal{G}_{i}}=\vword_{(J_i,b)},
 \end{equation}
  for $0\leq i \leq s$, can be defined.
  In order to unify notation we call $(V, \mathcal{G}, E)$ a "basic object", and a sequence of transformations  (\ref{sectransp2}) will be denoted by:
 \begin{equation}\label{sectrap7}
\begin{array}{cccccccc}
 & (V,\mathcal{G}_{(J,b)},E) & \stackrel{\pi_1}{\longleftarrow} &( V_{1},\mathcal{G}_{(J_1,b)},E_1) & \stackrel{\pi_2}{\longleftarrow}
 &\ldots& \stackrel{\pi_k}{\longleftarrow}&(V_s, \mathcal{G}_{(J_s,b)}, E_s).\\
 \end{array}
\end{equation}

  If each center, say $Y_i$, of $\pi_i$ is such that
  $Y_i \subset \vMax \vword_{\mathcal{G}_{i}}(=\vMax \vword_{(J_i,b)}) $,
then $$ \max \vword_{\mathcal{G}_{}} \geq \max \vword_{\mathcal{G}_{1}}\geq \dots \geq \max
\vword_{\mathcal{G}_{s}}$$
}
\end{parrafo}
\begin{parrafo}\label{qtts}{\em
 A Rees algebra $\mathcal{G}$ can be identified, up to integral closure, with two Rees algebras, say
$\mathcal{G}_{(J,b)}$ and $\mathcal{G}_{(J',b')}$, if and only if $(J,b)$ and $(J',b')$ are integrally equivalent. In particular satellite functions are well defined for $\mathcal{G}$ .

We say that  $\mathcal{G}_s$ is {\em monomial} if $\max \vword_{\mathcal{G}_s}=0$. This amounts to saying that $(J_s,b)$ is monomial (and that $(J'_s,b')$ is monomial).

 It is clear that also the second coordinate of the satellite function in \ref{poift}, and hence the inductive function itself,  extends naturally to the case of Rees algebras. In particular we say that (\ref{sectransp2}) is $t$-permissible if and only if  (\ref{resol341pp}) is $t$-permissible; namely if
  $$Y_i \subset \vMax t_i(\subset \vMax \vword_{\mathcal{G}_{i}}(=\vMax \vword_{(J_i,b)})),  $$
  in which case
  $ \max t\geq \max t_1\geq \dots \geq \max t_s$
  }
\end{parrafo}

 \begin{definition}
We say that $\mathcal{G}$ is {\em simple} (or that $(V,\mathcal{G}, E=\emptyset)$ is simple) if
$ord_{\mathcal{G}}(x)=1$ for all $x\in \Sing(\mathcal{G})$.
\end{definition}
 \begin{parrafo}\label{fgyz}{\em

Assume that a $t$-permissible sequence of monoidal transformations is defined, say:
 \begin{equation}\label{swt4ap7}
\begin{array}{cccccccc}
 & (V,\mathcal{G}_{},E) & \stackrel{\pi_1}{\longleftarrow} &( V_{1},\mathcal{G}_{1},E_1) & \stackrel{\pi_2}{\longleftarrow}
 &\ldots& \stackrel{\pi_k}{\longleftarrow}&(V_s, \mathcal{G}_{s}, E_s).\\
 \end{array}
\end{equation}

i) If $\max  \vword_s>0$, set $r$ as the smallest index so that
$\max  t_r=\max  t_{r+1}=\cdots = \max  t_s.$

ii) If $\max  \vword_s=0$, set $r$ as the smallest index so that
$\max  \vword_r=0$.

The following proposition is simply a reformulation of Prop \ref{propdegbo}.
}
\end{parrafo}

\begin{proposition} \label{Qropdegbo} Let (\ref{swt4ap7}) be a $t$-permissible sequence of monoidal transformations, and assume  that $\max  \vword_s>0$.  Fix $r$ as above. There is a simple Rees algebra $\mathcal{G}'$ at $V_r$, so that $(V_r, \mathcal{G}',E'_r=\emptyset)$ has the following property:

 Any sequence of transformations of $(V_r,\mathcal{G}',E'_r=\emptyset)$, say
\begin{equation}
\label{rezolpq6}
\begin{array}{cccccccc}
(V_r,\mathcal{G}',\emptyset)& \leftarrow &(V'_{r+1},(\mathcal{G}')_1,E'_{r+1}) & \leftarrow &\cdots
&\leftarrow & (V'_S,(\mathcal{G}')_S,E'_S),&
\\
\end{array}
\end{equation}
induces a $t$-permissible sequence of transformation of the basic object $(V_r,\mathcal{G}_r,E_r)$, say:
\begin{equation}
\label{qso76}
\begin{array}{cccccccc}
(V_r,\mathcal{G}_r,E_r)& \leftarrow &(V'_{r+1}, \mathcal{G}_{r+1},E_{r+1}) & \leftarrow &\cdots
&\leftarrow & (V'_S,\mathcal{G}_{S},E_S),&
\\
\end{array}
\end{equation}
by blowing up at the same centers, so that condition a), b), and c) of \ref{propdegbo} (suitable adapted) hold.

\end{proposition}
\proof Take  $(V_r,(J'_r,b'),\emptyset)$ as in Prop \ref{propdegbo}, and set
$\mathcal{G}'=\mathcal{G}_{(J',b')}.$

\begin{parrafo}{\bf  Rees algebras and pull-backs} {\em The analogy between the notions of Rees algebras and that of basic objects is also preserved by pull-backs. Pull-backs, defined for basic objects in (\ref{eqdrtbo1}),  can be reformulated as:
\begin{equation}\label{eqdrzo1}
 (V,\mathcal{G},E) \stackrel{\pi}{\longleftarrow}  (U,\mathcal{G}_U,E_U)
\end{equation}
which we call again a pull-back, which essentially is a restriction to an open set, or a restriction followed by multiplication by an affine space. We reformulate Definition \ref{lsotr}:
\begin{definition} \label{lsot9r} A {\em local sequence of transformations} of a basic object
$(V,\mathcal{G}, E)$ is
\begin{equation}\label{ABCranso}
(V,\mathcal{G},E) \longleftarrow (V'_1,\mathcal{G},E_1)\longleftarrow \cdots
\longleftarrow (V'_s,\mathcal{G},E_s),
\end{equation}
where $(V,\mathcal{G},E) \longleftarrow (V'_1,\mathcal{G}_1,E_1)$, and each
$(V'_i,\mathcal{G}_i,E_i) \longleftarrow (V'_{i+1},\mathcal{G}_{i+1},E_{i+1})$,
 is a pull-back, or a pull-back followed by a usual transform of a basic object (blowing up a smooth center $Y_i \subset \Sing(\mathcal{G}_i$ with normal crossings with hypersurfaces in $E_i$).
\end{definition}
The formulation of Def \ref{deeqft} in the context of Rees algebras is:
\begin{definition}\label{de7qft}
Two Rees algebras $\mathcal{G}^{(i)}$ , $i=1,2$, or two basic objects $(V,  \mathcal{G}^{(i)},E)$ , $i=1,2$, are said to be {\em weakly equivalent} if:
$\Sing(\mathcal{G}^{(1)})=\Sing(\mathcal{G}^{(2)}),$ and if
\begin{equation*}
(V',\mathcal{G}^{(i)},E') \longleftarrow (V'_1,\mathcal{G}^{(i)}_1,E'_1)\longleftarrow \cdots
\longleftarrow (V'_s,\mathcal{G}^{(i)}_s,E_s),
\end{equation*}
is a local sequence of transformations of one of them, then it  also defines a local sequence of transformation of the other, and $\Sing(\mathcal{G}_j^{(1)})=\Sing(\mathcal{G}_j^{(2)})$, $0\leq j \leq s$.
\end{definition}
The following is essentially a corollary of Lemma \ref{th41}.
\begin{theorem} \label{hggt} Let $\mathcal{G}$ be a Rees algebra on $V$, then $\mathcal{G}$ and $G(\mathcal{G})$
(or the basic objects $(V,  \mathcal{G},E)$  and $(V,  G(\mathcal{G}),E)$) are weakly equivalent.
\end{theorem}
Hironaka's Finite Presentation Theorem applies for Diff-algebras:
\begin{theorem} \label{hfpt} If two Diff-algebras on $V$, say  $\mathcal{G}$ and $\mathcal{G}' $  are weakly equivalent ( or, say, if $(V,  \mathcal{G},E=\emptyset)$  and $(V,  \mathcal{G}',E=\emptyset)$ are weakly equivalent), then
$\mathcal{G}$ and $\mathcal{G}' $ are integrally equivalent .
\end{theorem}
\begin{corollary} Let $\mathcal{G}$ and $\mathcal{G}'$ be two weakly equivalent Rees algebras on $V$,  then $G(\mathcal{G})$ and $G(\mathcal{G}')$
are integrally equivalent.

\end{corollary}
}
\end{parrafo}

\begin{parrafo}{\em

The main formula  (\ref{eqdjtbo6}) (of compatibility with pull-backs) is now expressed as:
\begin{equation}\label{eqdjtbo7}
 ord^{}_{\mathcal{G}_U}(x)=ord^{}_{\mathcal{G}}(\pi(x)).
 \end{equation}
This will allow us to extend the satellite functions to the case of  {\em local  transformations} (see Def \ref{lsotr}). If now
 \begin{equation}\label{sectp97}
\begin{array}{cccccccc}
 & (V',\mathcal{G}_{},E) & \stackrel{\pi_1}{\longleftarrow} &( V'_{1},\mathcal{G}_{1},E_1) & \stackrel{\pi_2}{\longleftarrow}
 &\ldots& \stackrel{\pi_k}{\longleftarrow}&(V'_s, \mathcal{G}_{s}, E_s).\\
 \end{array}
\end{equation}
is a local sequence of transformations, where each
\begin{equation}\label{eq9ff}
(V'_i,\mathcal{G}_{i},E_i) \longleftarrow (V'_{i+1},\mathcal{G}_{i+1},E_{i+1})
\end{equation}
is the blow up at a center $Y_i\subset \vMax \vword_i$ followed by a pull-back, then
$$ \max \vword \geq \max \vword_1\geq \dots \geq \max
\vword_s.$$
A similar argument applies  for the satellite function $t$ : if the previous condition holds, and $\max \vword_s>0$, then $t_i$ can be defined for $0\leq i \leq s$. Furthermore, if each "local transformation"  (\ref{eqaff}) is the blow up at a center $Y_i\subset \vMax t_i$ followed by a pull-back, then
$$ \max t \geq \max t_1\geq \dots \geq \max
t_s.$$
In this case we shall say that (\ref{sectp97}) is a $t$-permissible sequence of "local transformations".

Part A) of the following Proposition state properties of the function $t$ which are stronger and imply  Proposition \ref{Qropdegbo}. Parte B) will lead us a precise answer to assertion 2) in \ref{comnt}.

}
\end{parrafo}

\begin{proposition} \label{prxybo1} Let  (\ref{sectp97}) be a  $t$-permissible sequence of local transformations, assume that  $\max  \vword_s>0$. Let  $r$ be the smallest index so that $\max t_r=\max t_s$.

A) There is a simple Rees algebra $\mathcal{G}'$ at $V_r$ (or say $(V_r, \mathcal{G}',E'_r=\emptyset)$ ), with the following property: Any {\em local sequence of transformations} of $(V_r, \mathcal{G}',E'_r=\emptyset)$, say
\begin{equation}
\label{reesolll6}
\begin{array}{cccccccc}
(V_r, \mathcal{G}',E'_r=\emptyset)& \leftarrow &(V'_{r+1},  \mathcal{G}'_1,E'_{r+1}) & \leftarrow &\cdots
&\leftarrow & (V'_S,  \mathcal{G}'_S ,E'_S),&
\\
\end{array}
\end{equation}
induces a $t$-permissible local sequence of transformation of the basic object $(V_r, \mathcal{G}_r,E_r)$, say:
\begin{equation}
\label{reesollll6}
\begin{array}{cccccccc}
(V_r, \mathcal{G}_r,E_r)& \leftarrow &(V'_{r+1}, \mathcal{G}_{r+1},E_{r+1}) & \leftarrow &\cdots
&\leftarrow & (V'_{r+S},\mathcal{G}_{r+S},E_{r+S}),&
\\
\end{array}
\end{equation}
with the following condition on the functions $t_j$ defined for this last sequence (\ref{reesollll6}):

a) $\vMax t_{r+k}=\Sing( \mathcal{G}'_k)$, $1\leq k \leq S-1$.

b) $\max  t_r=\max  t_{r+1}=\cdots = \max  t_{r+S-1} \geq \max  t_{r+S}.$

c) $\max  t_{r+S-1} = \max  t_{r+S}$ iff  $\Sing(\mathcal{G}'_S) \neq \emptyset$ , in which case $\vMax t_{r+S}=\Sing(\mathcal{G}'_S)$.

B) If a Rees algebra $\mathcal{G}''$ at $V_r$  also fulfills A), then $\mathcal{G}''$ and $\mathcal{G}'$ are weakly equivalent.
\end{proposition}
\proof
Part A) is a reformulation of Proposition \ref{prdegbo1} . Part B) follows from A) and the definition of weak equivalence.

\begin{parrafo}\label{parKol}{\em

Part B) of the previous result, together with Theorem \ref{hggt} say that we may take $\mathcal{G}'$ to be a (simple) Diff-algebra; furthermore, Hironaka's Theorem \ref{hfpt} says that up to integral closure, there is a unique Diff-algebra
 $\mathcal{G}'$  which fulfills 2').

In the previous Proposition we could have taken  the Rees algebra $\mathcal{G}'$ to be of the form $\mathcal{G}_{(J'_r,b')}$, for  $(J'_r,b')$ as in \ref{ammr}. The Diff-algebra $G(\mathcal{G}_{(J'_r,b')})$
also fulfills the property; and it is unique with this property up to integral closure.

Suppose now that we know how to resolve simple basic objects (or simple Rees algebras)  in an constructive way, so that two Rees algebras with the same integral closure undergo the same resolution.
In this case the discussion in Remark  \ref{rrr43} says that, given a basic object,  the inductive function $t$ defines a {\em unique sequence of monoidal transformations} which make the basic object monomial.

Set  $G(\mathcal{G}_{(J'_r,b')})=\bigoplus_{k\geq 0}I_rW^r$. It is integral over
a Rees ring of an ideal, say $\calo_V[I_NW^N]$ for suitable $N$
(see \ref{rkK1}). These ideals $I_N$ are called {\em tuned ideals}
in \cite{kollar} (see page 45). We may replace $(J'_r,b')$  in Prop \ref{propdegbo}), by "tuned couple"
$(I_N,N)$, and if two tuned couples fulfill  the conditions of $(J'_r,b')$  they must idealistic equivalent.
This answers 2) in \ref{comnt}.

}
\end{parrafo}


\section{Projection of differential algebras and elimination.}
\label{sec5}

As was indicated in \ref{parKol}, Proposition \ref{prxybo1} ensures that if we know how to resolve simple basic objects, then a sequence of monoidal transformations can be defined over a basic object, so as to bring it to a simplified form (to the monomial case).  It also indicates some form of {\em uniqueness} in such procedure, a property which must hold in any constructive or algorithmic resolution. In this section we generalize that Proposition in Proposition \ref{prxyb81}, and we make use of the notion of elimination algebras introduced in \cite{VV4}, and generalized in \cite{BV3} . This last result, together with Theorem \ref{hfpt}, will lead us to the upper semi-continuous functions that stratify the singular locus into locally closed smooth sets in Theorem \ref{thbv}.

\begin{parrafo}\label{slsdh}{\em
Let $\mathcal{G}=\oplus I_sT^s$  be a Rees algebra on the smooth scheme $V$. Fix a closed point $x\in \Sing(\mathcal{G})$, with residue field $k'$, and a regular system of parameters
$\{x_1,\dots,x_n\}$ at ${\calo}_{V,x}$. $gr_{M_x}( {\calo}_{V,x})$,
is a polynomial ring, say $k'[X_1,\dots , X_n]$, where $X_i$ denotes the initial form of $x_i$, and  $Spec(gr_{M_x}( {\calo}_{V,x}))=\mathbb{T}_{V,x}$ (tangent space).
There is, on the one hand, a Taylor morphism, say:
$Tay: {\calo}_{V,x}\to {\calo}_{V,x}[[T_1,\dots ,
T_n]]$ that map $x_i$ to $x_i+T_i$, $1\leq i \leq n$ (see \ref{ffbrd}); on the other hand there is
a Taylor morphism, say:
$$\overline{Tay}:  k'[X_1,\dots , X_n]\to  k'[X_1,\dots , X_n][T_1,\dots ,
T_n]$$ that map $X_i$ to $X_i+T_i$, $1\leq i \leq n$.

Both are closely related, and some important invariants in desingularization arise from the link among them. In both cases we define, for each multi-index
$\alpha \in
\mathbb{N}^n$, operators, say
$\Gamma^{\alpha}$ and $\overline{\Gamma}^{\alpha}$, so that
$Tay(f)= \sum_{\alpha \in
\mathbb{N}^n} \Gamma^{\alpha}(f)T^{\alpha}$, and $\overline{Tay}(F)=\sum_{\alpha \in
\mathbb{N}^n} \overline{\Gamma}^{\alpha}(F)T^{\alpha}$.

Note here that if $F$ is an homogeneous polynomial of degree $N$,  and if   $|\alpha| \leq N$,
then $\overline{\Gamma}^{\alpha}(F)$ is either zero or homogeneous of degree $N-|\alpha|$.

As in \ref{trnvr1}, we attach an homogeneous ideal to $\mathcal{G}$ at $x$, say
$In_x(\mathcal{G})$, included in
 $gr_{M_x}( {\calo}_{V,x})$ ; namely that generated by the class of $I_s$ at the quotient $M_x^s/M_x^{s+1}$, for all $s$.  This homogeneous ideal defines a cone, say $C_{\mathcal{G}}$, at $\mathbb{T}_{V,x}$. Recall that there is a biggest  subspace, say  $L_{\mathcal{G}}$,  included and acting by translations on $C_{\mathcal{G}}$ (see \ref{trnvr1}). Hironaka defines $\tau_{\mathcal{G}}(x)$ (the $\tau$-invariant at the point) to be the codimension of the subspace $L_{\mathcal{G}}$ in $\mathbb{T}_{V,x}$.

 Recall also that $\Sing(\mathcal{G})= \Sing(G(\mathcal{G}))$. The relation among the two Taylor morphisms discussed above show how the two homogeneous ideals at $x$, attached to $G(\mathcal{G})$ and to $\mathcal{G}$ respectively, are related (namely $In_x(\mathcal{G})$ and $In_x(G(\mathcal{G}))$) : If $C_\mathcal{G}$ is the cone attached to $\mathcal{G}$, then the cone attached to $G(\mathcal{G})$ is the linear subspace  $L_{\mathcal{G}}$.

 In fact, the  graded ideal $In_x(G(\mathcal{G}))$
  is the smallest homogeneous extension of  $In_x(\mathcal{G})$ , closed by the action of the differential operators $\overline{\Gamma}^{\alpha}$; namely, with the property that if $F$ is an homogeneous polynomial of degree $N$ in the ideal, and if $|\alpha|< N$, then also $\overline{\Gamma}^{\alpha}(F)$ is in the  ideal. This homogeneous ideal defines the subspaces $L_{\mathcal{G}}$,  included in $C_{\mathcal{G}}$, with the properties stated in \ref{trnvr1}. Note, in particular, that $\mathcal{G}$ and $G(\mathcal{G})$ have the same $\tau$-invariant at all singular point (see \cite{Hironaka70} , \cite{Hironaka77}, \cite{Oda1973} , \cite{Oda1983} , \cite{Oda1987} ).

}
\end{parrafo}

Definition \ref{3def1} has a natural formulation in the relative context, namely when $\beta: V\to V'$ is a smooth morphism.
\begin{definition} \label{3def11} We say that a Rees algebra $\bigoplus I_nW^n$, on a smooth scheme $V$,  is a Diff-algebra relative to $V'$, if: i) $I_n\supset I_{n+1}$. ii) $Diff^{(r)}(I_n)\subset I_{n-r}$ for each $n$, and $0\leq r
\leq n$, where $Diff^{(r)}$ denotes the sheaf of differentials relative to $\beta: V\to V'$.
\end{definition}
The smooth morphism $\beta: V\to V'$ defines, at each point $x\in V$, a linear map at tangent spaces: $d\beta_x: \mathbb{T}_{V,x} \to
\mathbb{T}_{V',\beta(x)}$; and the kernel, say $ker  \ d\beta_x$  is a linear subspace of $\mathbb{T}_{V,x} $.
\begin{definition} \label{3def12} Fix a Rees algebra $\mathcal{G}=\bigoplus I_nW^n$ on $V$,  and a closed point
$x\in\Sing(\mathcal{G})$. We say that $\beta: V\to V'$ is {\em transversal} to $\mathcal{G}$ at the point, if
the subspaces $L_{\mathcal{G}}$ and $ker  \ d\beta_x$ are in general position in $\mathbb{T}_{V,x} $.
\end{definition}

\begin{definition} Fix an integer $e\geq 0$. We say that a Rees algebra $\mathcal{G}$ on a $d$-dimensional scheme $V$ is of {\em codimensional type} $e$, or say $\tau_{\mathcal{G}}\geq e$,  if
$\tau_{\mathcal{G}}(x)\geq e$ for all $x\in \Sing(\mathcal{G})$.

\end{definition}
\begin{remark}\label{eemio}
 If $\tau_{\mathcal{G}}\geq e$, the codimension of the closed set $\Sing(\mathcal{G})$ in $V$ is at least $e$. The components of codimension $e$ are open and closed in $\Sing(\mathcal{G})$, and they are regular. Furthermore, if $\Sing(\mathcal{G})$ is of pure codimension $e$, the monoidal transform at such center is a resolution of $\mathcal{G}$. This is, essentially, how centers are chosen in constructive resolutions, and also the reason why strata are smooth in the stratifications defined by the algorithm of desingularization.
 The codimensional type $e$ can be at most $d$ (dimension of $V$), and if  $\mathcal{G}$ is of codimensional type $d$, then $\Sing(\mathcal{G})$ is a finite set of points and a resolution of $\mathcal{G}$ is defined by the blow up at those points.

\end{remark}

\begin{parrafo}\label{ekli}{\bf Elimination algebras.} {\em  Assume that:

i) $\mathcal{G}$ is a Rees algebra on $V$ of codimensional type $\geq e$; and that $\Sing(\mathcal{G})$ has no component of codimension $e$.

ii)  $V'$ is smooth, dim$V$-dim$V'=e$ (i.e. dim$V'=d-e$), and $\beta: V\to V'$ is smooth, and transversal to $\mathcal{G}$ locally at every closed point $x\in \Sing(\mathcal{G})$.

iii) $\mathcal{G}$ is a Diff-algebra relative to $\beta: V\to V'$ .

If these conditions hold, then a Rees algebra, say $\mathcal{G}_{\beta}^{(e)}$, will be defined at the smooth scheme $V'$. $\mathcal{G}_{\beta}^{(e)}$ is called the {\em elimination algebra} (see Def 4.10, \cite{VV4} for the case $e=1$, and  \cite{BV3} for the general case).

Recall that the linear subspace $L_\mathcal{G}$, attached to  $\mathcal{G}$ at $\mathbb{T}_{V,x}$, is the same as that attached to $G(\mathcal{G})$.
So i) and ii) hold for $\mathcal{G}$ iff they hold for $G(\mathcal{G})$. On the other hand iii) will hold for
$G(\mathcal{G})$ with independence of $\beta: V\to V'$ . The local condition at $x$ in ii) is that $L_{\mathcal{G}}\cap ker  \ d\beta_x=0$.  In other words, let $\mathcal{G}$ be a Rees algebra of codimensional type $\geq e$ that fulfills  condition i), then elimination algebras will be defined for $G(\mathcal{G})$ locally at any point of $\Sing(G(\mathcal{G}))=\Sing(\mathcal{G})$, for
any smooth map as in ii).
}
\end{parrafo}
\begin{parrafo}\label{ekliA}{\bf Elimination algebras and local transformations.}{\em The following properties ensure the compatibility of elimination algebras with monoidal transformations and with pull-backs:

Ai) If $Y$ is smooth and included in $\Sing(\mathcal{G})$, then $\beta(Y)$ is smooth in $V'$ and included in   $\Sing(\mathcal{G}_{\beta}^{(e)})$.

Aii) Let $V\leftarrow V_1$ be the blow up at $Y$, and let  $\mathcal{G}_1$ be the transform of $\mathcal{G}$. Set $V'\leftarrow V'_1$ to be the blow up at $\beta(Y)$, and  let $(\mathcal{G}_{\beta}^{(e)})'_1$ be the transform of $\mathcal{G}_{\beta}^{(e)}$.

 There is a smooth map $\beta_1: V_1\to V_1'$ so that i) , ii) , and iii) hold at an open neighborhood of
$ \Sing(\mathcal{G}_1)$, and the elimination algebra of $\mathcal{G}_1$ is $(\mathcal{G}_{\beta}^{(e)})'_1$.

B) We consider conditions i) , ii) ,and iii), for  a Rees algebra  $\mathcal{G}$ on a smooth scheme $V$, the notion of monoidal transform leads to the consideration of basic objects, say
$(V,\mathcal{G},E)$.

 Whenever  $(V,\mathcal{G},E) \stackrel{\sigma}{\longleftarrow}  (U,\mathcal{G}_U,E_U)$ is a pull-back (see \ref{eqdrzo1}) the local conditions i),ii), and iii), have a natural lifting to
$ (U,\mathcal{G}_U,E_U)$, and the elimination algebra of the pull-back $\mathcal{G}_U$ is the pull-back of the elimination algebra.
}
\end{parrafo}
\begin{parrafo}\label{ekliB}{\bf Elimination and singular loci.} {\em  Set $\mathcal{G}=\bigoplus I_nW^n$ and $\beta: V\to V'$  so that the three conditions in \ref{ekli} hold. Set $\mathcal{G}_{\beta}^{(e)}=\bigoplus L_nW^n \subset \calo_{V'}[W].$ The following conditions hold:

1) $\mathcal{G}_{\beta}^{(e)}(=\bigoplus L_nW^n)\subset \mathcal{G}(=\bigoplus I_nW^n)$ via the inclusion
$\calo_{V'}[W]\subset \calo_{V}[W]$ defined by $\beta$. We also denote this  by
$\beta^*(\mathcal{G}_{\beta}^{(e)})\subset \mathcal{G}$, where $\beta^*( \calo_{V'})\subset  \calo_{V}$ is the inclusion defined by $\beta$.

2) $\beta(\Sing(\mathcal{G}))\subset \Sing(\mathcal{G}_{\beta}^{(e)})$
and the induced map $\beta:  \Sing(\mathcal{G})\to \beta(\Sing(\mathcal{G}))$ is a bijection.

}
\end{parrafo}
\begin{parrafo}{\em

$ \mathcal{G}_{\beta}^{(e)}$ is a Rees-algebra on the smooth scheme $V'$, and there is a function $ord_{\mathcal{G}_{\beta}^{(e)}}$  defined on $\Sing(\mathcal{G}_{\beta}^{(e)})$ (see \ref{agre1}). We now define a function, say $ord_{\beta}^{(e)}$, on the closed set $\Sing(\mathcal{G})$ , as the restriction of $ord_{\mathcal{G}_{\beta}^{(e)}}$ to the subset $\beta(\Sing(\mathcal{G}))\subset \Sing(\mathcal{G}_{\beta}^{(e)})$, followed by the identification of
$ \Sing(\mathcal{G})$ with $\beta(\Sing(\mathcal{G}))$ in \ref{ekliB}, 2).

It turns out  that these function is independent of the choice of the smooth map $\beta$. This property will allow us to define a function $ord_{\mathcal{G}}^{(e)}(=ord_{\beta}^{(e)})$ on  $\Sing(\mathcal{G})$, for a  particular class of Rees-algebra.


Recall that elimination algebras are defined for Diff-algebra of codimensional type $\geq e$ locally at any singular point, for suitable smooth morphisms
(see \ref{ekli}). So a function $ord_{\mathcal{G}}^{(e)}$ could be defined for any Diff-algebra that fulfills
\ref{ekli}, i). Moreover, the properties in \ref{ekliA} show that elimination algebras are also defined  for monoidal transformations and for pull-backs of Diff-algebras. In particular functions $ord_{\mathcal{G}}^{(e)}$ could be defined also for monoidal transformations and for pull-backs of Diff-algebras.
 Furthermore, they are defined for a successive sequence of pull-backs and transformations of a Diff-algebra. So the class of Rees-algebras $\mathcal{G}$ for which  the functions $ord_{\mathcal{G}}^{(e)}$ are defined is closed under monoidal transformations and pull-backs.
The following generalizes Theorem 5.5 in \cite{VV4}.
}
\end{parrafo}

\begin{theorem} \cite{BV3} \label{tafbta}

a) Let $\mathcal{G}$ be a Rees algebra over a smooth scheme $V$ of dimension $d$, and
let $\beta: V \to V'$ and $\delta: V \to V''$ be smooth maps on smooth schemes $V'$ and $V''$ of the same dimension $d-e$. If the three conditions in \ref{ekli} hold for both smooth maps, then
$$ord_{\beta}^{(e)}=ord_{\delta}^{(e)}$$
as functions on  $\Sing(\mathcal{G})$ . Therefore a function $ord_{\mathcal{G}}^{(e)}\to \mathbb{Q}$
is well defined.

b) If $\mathcal{G}$ and $\mathcal{G}'$ have the same integral closure, and if the functions
$ord_{\mathcal{G}}^{(e)}$ and $ord_{\mathcal{G}}^{(e)}$ are defined, then $ord_{\mathcal{G}}^{(e)}=ord_{\mathcal{G}'}^{(e)}$ on $\Sing(\mathcal{G})=\Sing(\mathcal{G}')$.

c) If $ord_{\mathcal{G}}^{(e)}$ is defined for $\mathcal{G}$ and
$(V,\mathcal{G},E) \stackrel{\sigma}{\longleftarrow}  (U,\mathcal{G}_U,E_U)$ is a pull-back , then
$ord_{\mathcal{G}_U}^{(e)}$ is defined and
$ord_{\mathcal{G}_U}^{(e)}(x)=ord_{\mathcal{G}}^{(e)}(\sigma(x))$ for any
$x\in \Sing(\mathcal{G}_U)$.

\end{theorem}
\begin{parrafo}{\em  Theorem \ref{tafbta} enables us to define a function, say $ord_{\mathcal{G}}^{(e)}(=ord_{\beta}^{(e)})$ on  $\Sing(\mathcal{G})$.
As for the case of  $e=0$, note that every Rees algebra $\mathcal{G}$ on $V$ is of codimensional type $\geq 0$, and the conditions in \ref{ekli} hold for $\beta$ the identity map. Furthermore,  the previous function  $ord_{\mathcal{G}}^{(0)}$ is the usual function  $ord_{\mathcal{G}}$ (see \ref{agre1}).
In this case, we take any basic object $(V,\mathcal{G},E)$, and define the satellite functions in \ref{qtts}; which were defined entirely in terms of Hironaka's function $ord_{\mathcal{G}}$ (see also \ref{rk55}).

Recall also the notion of a $t$-permissible transformations in \ref{qtts}, and consider
$(V_r, \mathcal{G}',E'_r=\emptyset)$  as in Proposition \ref{Qropdegbo}. There $\mathcal{G}'$ is a simple Rees algebra, or equivalently of codimensional type $\geq 1$. The discussion in \ref{parKol} says that we may take $(V_r, \mathcal{G}',E'_r=\emptyset)$ to be a Diff-algebra, and in that case property A) of Proposition \ref{prxybo1} characterize this simple basic object up to integral closure (if  $(V_r, \mathcal{G}',E'_r=\emptyset)$ and $(V_r, \mathcal{G}'',E''_r=\emptyset)$  are Diff-algebras that fulfill A) they have the same integral closure).

In case $e\geq 1$, consider $\beta:V\to V'$, $ \mathcal{G}$ ( and $ \mathcal{G}_{\beta}^{(e)}$) so that the three conditions in \ref{ekli} hold. Take a basic objects $(V,\mathcal{G}_{},E) $ and  $(V',\mathcal{G}^{(e)},E) $ where we assume:

1) The hypersurfaces of $E$ in $V$ are the pull-back of the hypersurfaces of $E$ in $V'$ (pull-back via $\beta$).

2) That $\Sing(\mathcal{G})$ has no component of codimension $e$.

A sequence of permissible monoidal transformations, say
\begin{equation}\label{swt5zp7}
\begin{array}{cccccccc}
 & (V,\mathcal{G}_{},E) & \stackrel{\pi_1}{\longleftarrow} &( V_{1},\mathcal{G}_{1},E_1) & \stackrel{\pi_2}{\longleftarrow}
 &\ldots& \stackrel{\pi_k}{\longleftarrow}&(V_s, \mathcal{G}_{s}, E_s).\\
 \end{array}
\end{equation}
induces a sequence
\begin{equation}\label{swt5zp5}
\begin{array}{cccccccc}
 & (V',\mathcal{G}_{\beta}^{(e)},E) & \stackrel{\pi_1}{\longleftarrow} &( V'_{1},(\mathcal{G}_{\beta}^{(e)})_{1},E_1) & \stackrel{\pi_2}{\longleftarrow}
 &\ldots& \stackrel{\pi_k}{\longleftarrow}&(V'_s,( \mathcal{G}_{\beta}^{(e)})_{s}, E_s).\\
 \end{array}
\end{equation}
and for each index $i$ there is a smooth morphism $\beta_i: V_i\to V'_i$, so that the three conditions of \ref{ekli} hold, and $( \mathcal{G}_{\beta}^{(e)})_{i}$ is the elimination algebra of $\mathcal{G}^{}_{i}$ (\ref{ekliA}). So for each index $i$, there is an identification of $\Sing(\mathcal{G}^{}_{i})$ with $\beta_i(\Sing(\mathcal{G}^{}_{i}))$, and an inclusion $\beta_i(\Sing(\mathcal{G}^{}_{i}))\subset \Sing( \mathcal{G}^{(e)}_{i})$.
The function $ord^{(e)}_{\mathcal{G}_i}$, defined on $\Sing(\mathcal{G}^{}_{i})$, is by definition the restriction to  $\beta_i(\Sing(\mathcal{G}^{}_{i}))$ of the function $ord^{}_{\mathcal{G}^{(e)}_i}$.
In particular the satellite functions of the functions $ord^{}_{\mathcal{G}^{(e)}_i}$ give rise to satellite functions of  $ord^{(e)}_{\mathcal{G}_i}$, and we say that (\ref{swt5zp7}) is $t^{(e)}$-permissible if
(\ref{swt5zp5}) is $t$-permissible in the usual sense.

In the previous discussion we have considered only monoidal transformations. But the same holds for pull-backs, namely there is a compatibility of smooth maps and elimination algebras with pull-backs.
The notions of satellite functions and of $t^{(e)}$-permissibility extend to the case of local transformations.

}
\end{parrafo}

\begin{proposition} \label{prxyb81} With the assumptions and hypothesis stated above, let
\begin{equation}\label{swt5qp7}
\begin{array}{cccccccc}
 & (V,\mathcal{G}_{},E) & \stackrel{\pi_1}{\longleftarrow} &( V_{1},\mathcal{G}_{1},E_1) & \stackrel{\pi_2}{\longleftarrow}
 &\ldots& \stackrel{\pi_k}{\longleftarrow}&(V_s, \mathcal{G}_{s}, E_s).\\
 \end{array}
\end{equation}
be $t^{(e)}$-permissible local sequence of transformations, assume that  $\max  \vword^{(e)}_s>0$, and fix  $r$ as in Prop \ref{Qropdegbo} (smallest index so that $\max t_r^{(e)}=\max t_s^{(e)}$).

A) There is a Rees algebra $\mathcal{G}''$ at $V_r$ (or say $(V_r, \mathcal{G}'',E'_r=\emptyset)$ ),  of codimensional type $\geq e+1$ with the following property:  Any {\em local sequence of transformations} of $(V_r, \mathcal{G}'',E'_r=\emptyset)$, say
\begin{equation}
\label{reesolll6}
\begin{array}{cccccccc}
(V_r, \mathcal{G}'',E'_r=\emptyset)& \leftarrow &(V'_{r+1},  \mathcal{G}''_1,E'_{r+1}) & \leftarrow &\cdots
&\leftarrow & (V'_S,  \mathcal{G}''_S ,E'_S),&
\\
\end{array}
\end{equation}
induces a $t^{(e)}$-permissible local sequence of transformation of  $(V_r, \mathcal{G}_r,E_r)$, say:
\begin{equation}
\label{reesollll6}
\begin{array}{cccccccc}
(V_r, \mathcal{G}_r,E_r)& \leftarrow &(V'_{r+1}, \mathcal{G}_{r+1},E_{r+1}) & \leftarrow &\cdots
&\leftarrow & (V'_{r+S},\mathcal{G}_{r+S},E_{r+S}),&
\\
\end{array}
\end{equation}
with the following condition on the functions $t^{(e)}_j$ defined for this last sequence (\ref{reesollll6}):

a) $\vMax t^{(e)}_{r+k}=\Sing( \mathcal{G}''_k)$, $1\leq k \leq S-1$.

b) $\max  t^{(e)}_r=\max  t^{(d-e)}_{r+1}=\cdots = \max  t^{(e)}_{r+S-1} \geq \max  t^{(e)}_{r+S}.$

c) $\max  t^{(e)}_{r+S-1} = \max  t^{e)}_{r+S}$ if and only if  $\Sing(\mathcal{G}''_S) \neq \emptyset$ , in which case $\vMax\mathcal{G}''=\Sing(\mathcal{G}''_S)$.

B) If $\mathcal{G}'''$ on $V_r$ also fulfills A), then it is weakly equivalent with $\mathcal{G}''$.
\end{proposition}
\proof A) Recall that (\ref{swt5zp7})  induces a sequence (\ref{swt5zp5}), and that there is a natural identification  of $\Sing( \mathcal{G}_i)$ with a closed subset  $\beta_i(\Sing( \mathcal{G}_i))$ of $\Sing( \mathcal{G}^{(e)}_i)$, for $0\leq i \leq s$. The function $t$ is upper semi-continuos, so after replacing each $V'_i$ by a suitable open neighborhood of  $\Sing( \mathcal{G}_i)$ we may assume that (\ref{swt5zp5}) is local $t$-permissible, that
 $\max  \vword_s>0$, and that $r$ is the smallest index so that $\max t_r=\max t_s$. Proposition \ref{prxybo1} applies to  (\ref{swt5zp5}), so let  $(V'_r, \mathcal{G}',E'_r=\emptyset)$ have the property stated in A) of that Proposition. The smooth morphism $\beta_r: V_r \to V'_r$ defines a lifting, say $\beta_r^*( \mathcal{G}')$, which is a Rees algebra on $V_r$  (see \ref{ekliB}).  Set $\mathcal{G}_{\beta}''=\mathcal{G}_r \odot \beta^*_r( \mathcal{G}')$ (see \ref{amalg}).
Note that $\Sing(\mathcal{G}_{\beta}'')=\vMax t^{(d-e)}_r$.

Let $Y\subset \Sing(\mathcal{G}_{\beta}'')$ be smooth, and define a transformation, say  $V_r \leftarrow W$, and transforms, say $(\mathcal{G}_{\beta}'')_1$ of $\mathcal{G}_{\beta}''$, and say $(\mathcal{G}_r)_1$ of $\mathcal{G}_r$. Then $Y$ is permissible for $\mathcal{G}_r$, and transversality ensures that $\beta_r(Y)$ is smooth at $V'_r$, and included in $\Sing( \mathcal{G}')$. Let $V'_r \leftarrow W'$ be the monoidal transformation and let
$(\mathcal{G}')_1$ be the transform of $\mathcal{G}'$. Then there is a natural lifting of $\beta_r$, say
$\beta'_r: W \to W'$, and one can finally check that
$(\mathcal{G}_{\beta}'')_1=(\mathcal{G}_r)_1\odot \beta'^{*}_{r}( (\mathcal{G}')_1)$. So again
$\Sing((\mathcal{G}_{\beta}'')_1)=\vMax t^{(e)}\subset ( \mathcal{G}_r)_1$.  A similar argument applies for pull-backs, and for any local permissible  sequence of $\mathcal{G}_{\beta}''$. So we can set
$\mathcal{G}_{}''=\mathcal{G}_{\beta}''$.

B) follows by definition, just as B) in Proposition \ref{prxybo1}.

\begin{remark} \label{llrk62}

 The definition of $\mathcal{G}_{\beta}''$ in Part A) of the previous proof was done in terms of the sequence (\ref{swt5zp5}) and the smooth morphisms $\beta_i: V_i \to V'_i$. Replace (\ref{swt5zp5}) by, say
\begin{equation}\label{swt5zpp}
\begin{array}{cccccccc}
 & (V'',\mathcal{G}^{(e)}_{\delta},E) & \stackrel{\pi_1}{\longleftarrow}
 &( V''_{1},(\mathcal{G}^{(e)}_{\delta})_1,E_1) & \stackrel{\pi_2}{\longleftarrow}
 &\ldots& \stackrel{\pi_k}{\longleftarrow}&(V''_s, (\mathcal{G}^{(e)}_{\delta})_s, E_s).\\
 \end{array}
\end{equation}
where for each index $i$ there is a smooth morphism $\delta_i: V_i\to V''_i$, so that the three conditions of \ref{ekli} hold, and $( \mathcal{G}^{(e)}_{\delta})_i$ is the elimination algebra of $\mathcal{G}^{}_{i}$.
 It follows from B)  that the Rees algebras $\mathcal{G}_{\beta}''$ and $\mathcal{G}_{\delta}''$ (defined in the proof of A), are weakly equivalent. In particular $G(\mathcal{G}_{\beta}'')$ and $G(\mathcal{G}_{\delta}'')$ have the same integral closure (Theorem \ref{hfpt}).

 \end{remark}
 \begin{remark}\label{llmon} Assume, as before that (\ref{swt5zp7}) is $t^{(e)}$-permissible and that $\max \vword_s^{(e)}=0$. This means that $ \mathcal{G}^{(e)}_{s}$ is a monomial Rees algebra (\ref{qtts}), at least in an open neighborhood of the closed set
 $\beta_s(\Sing(\mathcal{G}^{}_{s})) (\subset \Sing( \mathcal{G}^{(e)}_{s}))$. As the functions are independent of the projections, the same argument applies when we replace  (\ref{swt5zp5}) by
 (\ref{swt5zpp}): namely that  $ \mathcal{G}'^{(e)}_{s}$ is a monomial Rees algebra in an neighborhood
 of  $\delta_s(\Sing(\mathcal{G}^{}_{s}))$.
 \end{remark}
 \begin{corollary} Assume that one can define a resolution for any Diff-algebra of codimensional type
 $\geq e+1$. Let $\mathcal{G}$ be of codimensional type  $\geq e$ on a smooth scheme $V$. Then a $t^{(e)}$-permissible sequence, now of monoidal transformations, say
 \begin{equation}\label{wwt5zp7}
\begin{array}{cccccccc}
 & (V,\mathcal{G}_{},E) & \stackrel{\pi_1}{\longleftarrow} &( V_{1},\mathcal{G}_{1},E_1) & \stackrel{\pi_2}{\longleftarrow}
 &\ldots& \stackrel{\pi_k}{\longleftarrow}&(V_s, \mathcal{G}_{s}, E_s).\\
 \end{array}
\end{equation}
can be defined so that $\Sing(\mathcal{G}_{s})=\emptyset$ or $\max \vword_s^{(e)}=0$.

 \end{corollary}
 \begin{parrafo}\label{llmonn}{\em  The discussion in Remark \ref{llmon} indicates that  $\max \vword_s^{(e)}=0$ when the elimination algebra is (locally) a monomial Rees algebra in the sense of (\ref{qtts}). We will say that an e-codimensional Rees-algebra $ \mathcal{G}_{s}$ is  {\em e-monomial }  when $\max \vword_s^{(e)}=0$. So the previous Corollary says that by decreasing induction on $e$, we can define (\ref{wwt5zp7}) to be either a resolution, or $ \mathcal{G}_{s}$ is an {\em e-monomial }  Rees algebra.}

 \end{parrafo}
\begin{parrafo}{\em
 Set $\mathbb{T}=\mathbb{Q}\times \mathbb{Z}\cup
\{{\infty}\}$ so that $\mathbb{Q}\times \mathbb{Z}$ is ordered lexicographically, and $\{{\infty}\}$ is the biggest element. And set $I_d=\mathbb{T}\times \mathbb{T}\times
\dots \times \mathbb{T}$ (d-times $\mathbb{T}$) with lexicographic order.
}
\end{parrafo}
\begin{theorem}(\cite{BV3})  \label{thbv} Let $\mathcal{G}$ be a Diff-algebra on a smooth
scheme $V$ of dimension $d$. There is an upper semi-continuous function, say
$\gamma_{\mathcal{G}}: \Sing(\mathcal{G})\to I_d$ so that:

i) The level sets of $\gamma_{\mathcal{G}}$ stratify
$\Sing(\mathcal{G})$ in smooth (locally closed) strata.

ii) Over fields of characteristic zero, the function coincides
with the desingularization function used in \cite{Villa89}.

\end{theorem}
The proof follows from the following result:
\begin{proposition} \label{prop617}  Given a basic object $(V,\mathcal{G}_{},E)$, there is a unique sequence
\begin{equation}\label{ww3rp7}
\begin{array}{cccccccc}
 & (V,\mathcal{G}_{},E) & \stackrel{\pi_0}{\longleftarrow} &( V_{1},\mathcal{G}_{1},E_1) & \stackrel{\pi_1}{\longleftarrow}
 &\ldots& \stackrel{\pi_{s-1}}{\longleftarrow}&(V_s, \mathcal{G}_{s}, E_s).\\
 \end{array}
\end{equation}
together with upper semi-continuous functions $\gamma_{i}: \Sing(\mathcal{G})_i\to I_d$, so that for each index $i<s$ $\pi_i$ is the blow-up at the smooth scheme $Y_i=\vMax \gamma_i$, and either (\ref{ww3rp7}) is a resolution of the basic object, or $\vMax \gamma_s$ is the singular locus of an e-monomial Rees algebra (see \ref{llmonn} ) for some integer $e>0$.

Furthermore, for each index $i<s$ there is a positive integer $e_i$ and an $e_i$-codimensional Diff-algebra $\mathcal{G}''_i$ on $V_i$, so that $Y_i=\vMax \gamma_i$ is the union of components  of $\Sing(\mathcal{G}''_i)$ of codimension $e_i$  in $V_i$ (see \ref{eemio}).

\end{proposition}
\proof The inductive function $t$ is upper semi-continuous (\ref{poift}), so if we fix $x\in \Sing(\mathcal{G})$, the value, say $\alpha_1=t(x)\in \mathbb{Q}\times \mathbb{Z}$ is the highest value in a neighborhood of $x$.  Proposition \ref{prxyb81} and Remark \ref{llrk62}  (or say Theorem \ref{hfpt} ) allow us to define a {\em unique} Diff-algebra attached to the value $\alpha_1$, say $ \mathcal{G}_{\alpha_1}$   of  codimensional type 1.

 Define $\gamma(x)=(\alpha_1, \infty, \dots ,\infty)\in I_d$ if $x$ is in a component of codimension 1 of  $\Sing(\mathcal{G}_{\alpha_1})$. If not, an upper semi-continuous function $t^{(1)}$ is defined along $\Sing(\mathcal{G}_{\alpha_1})$; so set $\alpha_2=t^{(1)}(x)\in \mathbb{Q}\times \mathbb{Z}$ , and now there is a Diff-algebra $ \mathcal{G}_{{\alpha_1, \alpha_2}}$   of  codimensional type 2 attached to the value.

 Define $\gamma(x)=(\alpha_1, \alpha_2,\infty, \dots ,\infty)\in I_d$ if $x$ is in a component of codimension 2 of  $\Sing(\mathcal{G}_{{\alpha_1, \alpha_2}})$. If not, an upper semi-continuous function $t^{(2)}$ is defined at $x$, and we argue as before. In this way we define $\gamma$ along $\Sing(\mathcal{G})$. This function is upper semi-continuous, and has the following property. Assume that
$(\alpha_1, \alpha_2,\dots, \alpha_r,\infty, \dots ,\infty)$, is the highest value, and fix
$r'\leq r$. The set of points where where the first $r'$ coordinates of the function take the value
$(\alpha_1, \alpha_2,\dots, \alpha_{r'})$ is $ \Sing(\mathcal{G}_{{\alpha_1, \alpha_2, \dots , \alpha_{r'}}})$ where $\mathcal{G}_{{\alpha_1, \alpha_2, \dots , \alpha_{r'}}}$ is of codimensional type $\geq r'$. And for $r'=r$, $ \Sing(\mathcal{G}_{{\alpha_1, \alpha_2, \dots , \alpha_{r}}})$ is of pure codimension $r$, and this will be our choice of center for $\pi_0$. So the function define the center. The sequence (\ref{ww3rp7})  can be defined with the following property which makes it unique.

Assume, by induction on $s$, that a sequence as (\ref{ww3rp7}) is defined together with the functions
$\gamma_{i}: \Sing(\mathcal{G})_i\to I_d$. Set $\max \gamma_s=(\alpha_1, \alpha_2,\dots, \alpha_r,\infty, \dots ,\infty)$ (highest value achieved by $\gamma_s$) and fix $r'\leq r$. Let $s'$ be the smallest index for which the first $r'$ coordinates of $\max \gamma_{s'}$ is $(\alpha_1, \alpha_2, \dots , \alpha_{r'})$.
There is a Diff-algebra, say $\mathcal{G}_{{\alpha_1, \alpha_2, \dots , \alpha_{r'}}}$  of codimensional type $\geq r'$ at $V_{s'}$, so that
$ \Sing(\mathcal{G}_{{\alpha_1, \alpha_2, \dots , \alpha_{r'}}})$ are the points where the first $r'$ coordinates of $\gamma_{s'}$ take this value. The same centers of  transformations $\pi_i$, $i=s', \dots ,s$ in (\ref{ww3rp7}) define
\begin{equation*}
\begin{array}{cccccccc}
 & (V_{s'},\mathcal{G}_{{\alpha_1, \alpha_2, \dots , \alpha_{r'}}},E'_{s'}=\emptyset) & \stackrel{\pi_{s'}}{\longleftarrow} &( V_{s'+1},(\mathcal{G}_{{\alpha_1, \alpha_2, \dots , \alpha_{r'}}})_1,E'_{s'+1}) &
 &\ldots& \stackrel{\pi_k}{\longleftarrow}&(V_s, (\mathcal{G}_{{\alpha_1, \alpha_2, \dots , \alpha_{r'}}})_{s-s'}, E'_s).\\
 \end{array}
\end{equation*}
and $ \Sing((\mathcal{G}_{{\alpha_1, \alpha_2, \dots , \alpha_{r'}}})_{s-s'})$ is the set of points where the first
$r'$ coordinates of $\gamma_s$ take the value $(\alpha_1, \alpha_2,\dots, \alpha_{r'})$. Furthermore, this sequence is $t^{(r')}$-permissible and the $r'+1$-coordinate of $\gamma_s$ along $ \Sing((\mathcal{G}_{{\alpha_1, \alpha_2, \dots , \alpha_{r'}}})_{s-s'})$ is defined in terms of the function $t^{(r')}$.

If, in the previous discussion,  we take $r'=r$, then $(\mathcal{G}_{{\alpha_1, \alpha_2, \dots , \alpha_{r}}})_{s-s'}$ is $r$-codimensional, its singular locus is $\vMax \gamma_s$, and (\ref{ww3rp7}) can be extended if the function $\vword^{(r)}$, defined in terms of
$(\mathcal{G}_{{\alpha_1, \alpha_2, \dots , \alpha_{r}}})_{s-s'}$  is not zero (i.e. if $(\mathcal{G}_{{\alpha_1, \alpha_2, \dots , \alpha_{r}}})_{s-s'}$ is not $r$-monomial in the sense of Remark \ref{llmonn}).

\end{document}